\documentclass{amsart}
\usepackage{amssymb}
\usepackage{graphicx}

\usepackage[colorlinks]{hyperref}
\newtheorem{thm}{Theorem}[section]
\newtheorem{cor}[thm]{Corollary}
\newtheorem{thm-def}[thm]{Theorem-Definition}
\newtheorem{lem}[thm]{Lemma}
\newtheorem{prop}[thm]{Proposition}
\theoremstyle{definition}
\newtheorem{defn}[thm]{Definition}
\theoremstyle{remark}
\newtheorem{rem}[thm]{Remark}
\numberwithin{equation}{section}

\newcommand{\Spec}{{\rm Spec}}

\newcommand{\tensor}{\otimes}

\newcommand{\vacB}{|0\rangle}
\newcommand{\vac}{\vacB}
\newcommand{\Cplx}{\mathbb C}
\newcommand{\ZZ}{\mathbb Z}
\newcommand{\LL}{\mathbb L}
\newcommand{\VV}{\mathbb V}
\newcommand{\MM}{\mathbb M}
\newcommand{\gr}{{\rm gr\ }}
\newcommand{\grG}{{\rm gr^{\mathcal{G}}\ }}

\newcommand{\grH}{{\rm gr^{\mathcal{H}}\ }}

\newcommand{\pole}{\Cplx}

\newcommand{\de}{\partial}
\newcommand{\oX}{{\mathcal{O}}_{X}}

\newcommand{\calG}{{\mathcal{G}}}

\newcommand{\calT}{{\mathcal{T}}}

\newcommand{\cA}{\mathcal{A}}

\newcommand{\cD}{{\mathcal{D}}}

\newcommand{\cF}{{\mathcal{F}}}
\newcommand{\cG}{{\mathcal{G}}}
\newcommand{\cH}{{\mathcal{H}}}
\newcommand{\cJ}{{\mathcal{J}}}

\newcommand{\cL}{{\mathcal{L}}}
\newcommand{\cM}{{\mathcal{M}}}

\newcommand{\cO}{{\mathcal{O}}}

\newcommand{\cT}{{\mathcal{T}}}

\newcommand{\cV}{{\mathcal{V}}}

\newcommand{\cZ}{{\mathcal{Z}}}

\newcommand{\zero}{{}_{(0)} }
\newcommand{\one}{{}_{(1)} }
\newcommand{\opm}{{}_{(-1)} }
\newcommand{\ops}[1]{{}_{({#1})} }

\newcommand{\End}{{\mathrm End\,}}

\newcommand{\id}{{\rm id\ }}

\newcommand{\Ind}{{\rm Ind\ }}

\newcommand{\spn}{{\rm span}}
\newcommand{\Der}{{\rm Der\ }}

\newcommand{\iso}{\simeq}

\newcommand{\la}{\lambda}

\newcommand{\fg}{\frak{g}}

\newcommand{\fn}{\frak{n}}
\newcommand{\fh}{\frak{h}}
\newcommand{\fz}{\frak{z}}

\newcommand{\prline}{\mathbb{P}^1}

\newcommand{\Symb}{{\rm Symb\  }}

\newcommand{\cdo}{{\cD}_X^{ch}}

\newcommand{\twcdo}{{\cD}_X^{ch, tw}}
\newcommand{\twcdoloc}{{\stackrel{\circ}{\cD}}_X^{ch, tw}}
\newcommand{\gtdo}{{\cD}_X^{tw}}

\newcommand{\btau}{\bar{\tau}}
\newcommand{\bxi}{\bar{\xi}}

\newcommand{\set}[1]{\left\{#1\right\}}
\newcommand{\Fields}{{Fields}}
\newcommand{\ghat}{{\hat{\fg}}}
\newcommand{\qf}[2]{\langle {#1}, {#2}  \rangle}

\newcommand{\uH}{\underline{\rm{H}}}

\begin{document}

\pagestyle{plain}

\title{Algebras of twisted chiral differential operators and affine localization of \boldmath{$\fg$}-modules}
\author{T.Arakawa\and D.Chebotarov \and F.Malikov}
\maketitle

\begin{abstract}
We propose a notion of algebra of {\it twisted} chiral differential
operators over  algebraic manifolds with vanishing 1st Pontrjagin
class. We show that such algebras possess families of modules
depending on infinitely many complex parameters, which we classify
in terms of the corresponding algebra of twisted differential
operators. If the underlying manifold is a flag manifold, our
construction recovers modules over an affine Lie algebra
parameterized by opers over the Langlands dual Lie algebra. The
spaces of global sections of ``smallest'' such modules are
irreducible $\ghat$-modules and all irreducible  $\fg$-integrable
$\ghat$-modules at the critical level arise in this way.
\end{abstract}

\section{Introduction}
\subsection{Algebras} Algebras of twisted differential operators (TDO) were proposed
by Bernstein and Beilinson \cite{BB1,BB2} as a tool to study
representation theory of simple complex Lie algebras. To give an
example, consider the projective line $\prline$ with an atlas
consisting of 2 copies of $\pole$ with coordinates $x$ and $y$ resp.
so that $y=1/x$. One has
\begin{equation}
\label{glueonp1}
\partial_y=-x^2\partial_x.
\end{equation}
This defines the tangent sheaf $\cT_{\prline}$; $\cT_{\prline}$ is a
Lie algebroid and its universal enveloping algebra is the algebra of
differential operators $\cD_{\prline}$.

This construction is twisted by postulating the following transition
function
\begin{equation}
\label{glueonp1-tw}
\partial_y=-x^2\partial_x+\lambda x, \lambda\in\pole.
\end{equation}
The result is the algebra of twisted differential operators
$\cD^{\lambda}_{\prline}$. It is isomorphic to $\cD_{\prline}$
locally, but not globally; for example, if $\lambda$ is an integer,
then $\cD^{\lambda}_{\prline}$ is the algebra of differential
operators acting on the sheaf $\cO(\lambda)$.

Such algebras of locally trivial twisted differential operators can
be defined for an arbitrary smooth algebraic variety $X$; their
isomorphism classes are in 1-1 correspondence with
$H^1(X,\Omega_X^{1,cl})$. Thus for each $\lambda\in
H^1(X,\Omega_X^{1,cl})$, there is an algebra $\cD^{\lambda}_{X}$.

This construction can be further generalized to include algebras
that are not isomorphic to $\cD_X$ even locally. These are
classified by the hypercohomology group
$H^1(X,\Omega^1_X\rightarrow\Omega^{2,cl}_X)$, and we obtain a
$\cD^{\lambda}_X$ for each $\lambda\in
H^1(X,\Omega^1_X\rightarrow\Omega^{2,cl}_X)$.

Introduced in \cite{MSV,GMS} -- and in \cite{BD1} in the language of
chiral algebras -- are algebras of chiral differential operators,
CDO; these are sheaves of vertex  algebras of a certain type that
resemble algebras $\cD_X$ in some respects. A CDO over $X$ may or
may not exist; in fact, it exists if and only if $ch_2(\cT_X)\in
H^2(X,\Omega^2_X\rightarrow\Omega^{3,cl}_X)$ equals 0. If it does,
then the isomorphism classes of CDO-s over $X$ are a torsor over
$H^1(X,\Omega^2_X\rightarrow\Omega^{3,cl}_X)$ -- note that the
degree has jumped in comparison with the case of twisted
differential operators.

One can argue, therefore, that all CDO-s are twisted, because there
is no distinguished one and, worse still,  there may be none at all.
Nevertheless, it is the purpose of this paper to introduce a class
of {\it twisted chiral differential operators}, TCDO, so that all of
the above CDO-s will appear untwisted.

To give a flavor of the construction, let us return to the case of
$X=\prline$. $\prline$ carries a unique up to isomorphism CDO,
$\cD^{ch}_{\prline}$; it is defined by means of the following
`chiralization' of (\ref{glueonp1}):
\begin{equation}
\label{glueonp1-ch}
\partial_y=-x_{(-1)}x_{(-1)}\partial_x-2\partial(x),
\end{equation}
where we have let ourselves  use freely some of vertex algebra and
CDO notation; for example, $x$ and $\partial_x$ are fields
associated (in some sense) to the coordinate and  derivation so
denoted, and $\partial(x)$ means the canonical vertex algebra {\it
translation operator} applied to $x$.

Next, one would like to find a chiral version of
(\ref{glueonp1-tw}). Writing simply
$\partial_y=-x_{(-1)}x_{(-1)}\partial_x-2\partial(x)+\lambda x$,
$\lambda\in\pole$, is possible but uninteresting and ultimately
unhelpful. It appears that the right thing to do is to chiralize not
any of $\cD^{\lambda}_X$ but their universal version, $\cD^{tw}_X$.
In the case of $\prline$, this means to define $\cD^{tw}_{\prline}$
as a $\cO_{\prline}\otimes\pole[\lambda]$-module using the same
(\ref{glueonp1-tw}) with $\lambda$ not a number but a variable.

The chiral version of this is as follows: replace $\pole[\lambda]$
with
$H_{\prline}=\pole[\lambda,\partial(\lambda),\partial^2(\lambda),...]$,
the commutative vertex algebra of differential polynomials on
$\pole$, and then define an algebra of twisted chiral differential
operators, $\cD^{ch,tw}_{\prline}$, to be the sheaf of vertex
algebras locally isomorphic to $\cD^{ch}_{\prline}\otimes
H_{\prline}$ with the following transition functions
\begin{equation}
\label{glueonp1-ch-tw}
\partial_y=-x_{(-1)}x_{(-1)}\partial_x-2\partial(x)+\lambda_{(-1)}x.
\end{equation}

Similarly, we construct for an arbitrary compact smooth $X$ the
universal algebra of twisted differential operators, $\cD^{tw}_X$;
it is an algebra over
$\pole[H^1(X,\Omega^1_X\rightarrow\Omega^{2,cl}_X)]$ such that being
quotiented out by  the maximal ideal of  a point $\lambda\in
H^1(X,\Omega^1_X\rightarrow\Omega^{2,cl}_X)$ it gives
$\cD^{\lambda}_X$. We then chiralize this construction and obtain,
for each CDO $\cD^{ch}_X$, a twisted CDO $\cD^{ch,tw}_X$, a sheaf of
vertex algebras, which locally, but not globally,  looks like
$\cD^{ch,tw}_X\otimes H_X$, where $H_X$ is the algebra of
differential polynomials on
$H^1(X,\Omega^1_X\rightarrow\Omega^{2,cl}_X)$.

Apart from serving as a prototype, algebras of twisted differential
operators are directly linked  to algebras of twisted chiral
differential operators via the notion of the Zhu algebra \cite{Zhu},
and this is another topic of the present paper. Zhu attached to each
graded vertex algebra $V$ an associative algebra, $\cZ hu(V)$.  We
show that the sheaf associated to the presheaf $X\supset U\mapsto
\cZ hu(\cD^{ch,tw}_X(U))$ is precisely $\cD^{tw}_X$.

$\cZ hu(V)$ controls representation theory of $V$, the subject to
which we now turn.

\subsection{Modules}

Note that $\cD^{ch,tw}_X$ is not a deformation of $\cD^{ch}_X$, not
technically at least, but it has a rich representation theory. In
particular, it has families of modules that are indeed deformations
of those over $\cD^{ch}_X$, and this is why $\cD^{ch,tw}_X$ may be
of interest.

Zhu showed that under some restrictions, a $V$-module is the same as
a $\cZ hu(V)$-module. It follows easily that (under similar
restrictions) a $\cD^{ch}_X$-module is the same as a $\cD_X$-module,
a result that is a bit disheartening.

One of those restrictions is that a $V$-module be graded. In the
case of $\cD^{ch,tw}_X$ let us relax this by demanding that modules
be only filtered. Now note that $H_X$ belongs to the center of
$\cD^{ch,tw}_X$. Therefore, we can take any $\cD^{ch,tw}_X$-module,
for example one coming from a $\cD^{tw}_X$-module, and quotient it
out by a character of $H_X$.

It is easy to see that  a character of $H_X$ is  an element of
$H^1(X,\Omega^1_X\rightarrow\Omega^{2,cl}_X)((z))$. Among those a
special role is played by characters with regular singularities,
$\chi(z)=\chi_0/z+\chi_{-1}+\chi_{-2}z+\cdots$, $\chi_{j}\in
H^1(X,\Omega^1_X\rightarrow\Omega^{2,cl}_X)$.

Arguing along these lines we prove that

\bigskip
{\it A $\cD^{ch,tw}_X$-module with central character $\chi(z)$ is
the same as a $\cD^{\chi_0}_X$-module
if $\chi(z)$ has regular
singularity and zero otherwise.}

%$\chi_0=\text{res}_{z=0}\chi(z)$,
\bigskip
The content of this assertion is  not in the vanishing result, which is valid
only under some technical restrictions that we have skipped anyway,
but in the explicit construction of a variety of modules labeled by
characters $\chi(z)$. Here is one example that this construction
generalizes.

Let $X=G/B$ be a flag manifold. Then a $\cD^{\chi_0}_X$-module is
essentially the same as a $\fg$-module with central character
determined by $\chi_0$. Applying the above construction to the
contragredient Verma module over $\fg$, we obtain a sheaf whose
space of sections over the big cell is the Wakimoto module over
$\widehat{\fg}$ at the critical level quotiented out by the central
character $\chi(z)$ \cite{FF1,F1}, $\chi(z)$   being interpreted in
this case as an {\it oper} for the Langlands dual group. Of course,
this is a beginning of the representation-theoretic input to the
Beilinson-Drinfeld construction of Hecke eigensheaves on $Bun_{G}$,
\cite{BD2}, also \cite{F2,F3}.  Therefore, what we are doing can be
thought of as providing ``operatic'' parameters in the case of an
arbitrary manifold; and indeed, the spectacular work by Feigin and
Frenkel served as a major source of inspiration for us.

Furthermore, we prove  that

\bigskip
{\it if the $\fg$-module we start with is simple and
finite-dimensional, then the space of global sections of the
corresponding $\cD^{ch,tw}_{G/B}$-module is irreducible and
isomorphic to the Weyl module over $\widehat{\fg}$ at the
critical level quotiented out by the central character.}

\bigskip
 The irreducibility of Weyl modules  at the
critical level quotiented out by the central character is a result
of Frenkel and Gaitsgory, which was anticipated in [FG2] and proved
in [FG3]. Our analysis of the spaces of global sections heavily
relies on techniques and results \cite{FG1,FG2,FG3}.

Let us see how this (and a bit more) comes about in the case of
$X=\prline$.

This case is described by explicit formulas (\ref{glueonp1-tw}) and
(\ref{glueonp1-ch-tw}). If we let in (\ref{glueonp1-tw})
$\lambda=n\in\ZZ$, then the ``smallest'' $\cD^{n}_{\prline}$-module
is $\cO(n)$. To make our life easier, let $\chi(z)=n/z$. This is the
case when the resulting $\cD^{ch,tw}_{\prline}$-module is actually
graded; denote it by $\cO(n)^{ch}$. Note that when $n=0$,
$\cO(0)^{ch}$ is precisely $\cD^{ch}_{\prline}$, and has been known
since \cite{MSV}.

We prove that

\bigskip
{\it(i) $H^0(\prline,\cO(n)^{ch})$ and $H^1(\prline,\cO(n)^{ch})$
are isomorphic to the irreducible $\widehat{sl}_2$-module at the
critical level with highest weight $n$ if $n\geq 0$;

(ii) $H^0(\prline,\cO(n)^{ch})$ and $H^1(\prline,\cO(n)^{ch})$ are
 isomorphic to
the irreducible $\widehat{sl}_2$-module at the critical level with
highest weight $-n-2$ if $n\leq -2$;

(iii) $H^0(\prline,\cO(n)^{ch})=H^1(\prline,\cO(n)^{ch})=0$  if
$n=-1$.}

\bigskip

This result is a direct generalization of \cite{MSV}, Theorem 5.7,
sect.5.8, and our construction verifies the proposals made in
\cite{MSV}, sect.5.15, one of the starting points of the present
work.

To conclude, one can say that the category of
$\cD^{ch,tw}_{G/B}$-modules appears to be a cross between the
Bernstein-Beilinson \cite{BB1,BB2} localization of $\fg$-modules to
the flag manifold and localization of $\widehat{\fg}$-modules at the
critical level to the semi-infinite flag manifold. We hope that this
point of view may prove useful.

\bigskip

{\it Acknowledgments.} We have benefited from discussions with P.Bressler, A.Beilinson,
D.Gaitsgory, V.Gorbounov, V.Schechtman. Part of this work was done when we were visiting the
Max-Planck-Institut f\"ur Mathematik in Bonn and Institut des Hautes \'Etudes Scientifiques in
Bures-sur-Yvette. We are grateful to these institutions for the superb working conditions. F.M.
was partially supported by an NSF grant. T. A. was partially supported by the JSPS Grant-in-Aid
for Scientific Research (B) No.\ 20340007.

\bigskip

\section{Preliminaries.}

We will recall the basic  notions of vertex algebra and describe
computational tools to be used in the sequel.

All vector spaces will be over  $\pole$. All spaces are even.
\subsection{Definitions.}
\label{Definitions.} Let $V$ be a vector space.

A {\em field} on $V$ is a formal series $$a(z) = \sum_{n\in \ZZ}
a_{(n)} z^{-n-1} \in ({\rm End} V)[[z, z^{-1}]]$$ such that for any
$v\in V$ one has $a_{(n)}v = 0$ for sufficiently large $n$.

Let $\Fields (V)$ denote the space of all fields on $V$.

A {\em vertex algebra}  is a vector space $V$ with the following
data:
\begin{itemize}
  \item a linear map $Y: V \to \Fields(V)$,  %  ({\rm End\  }V)[[z, z^{-1} ]]$,
    $V\ni a \mapsto a(z) = \sum_{n\in \ZZ} a_{(n)} z^{-n-1}$
  \item a vector $\vac\in V$, called {\em vacuum vector}
  \item a linear operator $\de: V \to V$, called {\em translation operator}
\end{itemize}
that satisfy the following axioms:
\begin{enumerate}
  \item (Translation Covariance)

 $ (\de a)(z) = \de_z a(z)$

  \item (Vacuum)

  $\vac(z) = \id$;

  $a(z)\vac \in V[z]$ and $a_{(-1)}\vac = a$

  \item (Borcherds identity)
   \begin{align}
   \label{Borcherds-identity}
   &\sum\limits_{j\geq 0} {m \choose j} (a\ops{n+j} b )\ops{m+k-j}\\
   =& \sum\limits_{j\geq 0} (-1)^{j} {n \choose j}\{ a\ops{m+n-j} b\ops{k+j} - (-1)^{n}b\ops{n+k-j} a\ops{m+j}
   \}\nonumber
   \end{align}
\end{enumerate}

% \end{defn}

A vertex algebra $V$ is {\em graded} if  $V = \oplus_{n\geq 0}V_n$
and for $a\in V_i$, $b\in V_j$ we have
$$a_{(k)}b \in V_{i+j - k -1}$$ for all $k\in \ZZ$. (We put $V_i  = 0$ for $i<0$.)

We say that a vector $v\in V_m$ has {\em conformal weight} $m$ and
write $\Delta_v = m$.

If $v\in V_m $ we denote $v_k  = v_{(k - m +1)}$, this is the
so-called conformal weight notation for operators. One has $$v_k V_m
\subset V_{m -k}.$$

 A {\em morphism} of vertex algebras is a map $f: V \to W$ that preserves vacuum and satisfies $f(v_{(n)}v') = f(v)_{(n)}f(v')$.

% \begin{defn}

A {\em module} over a vertex algebra $V$ is a vector space $M$
together with a map
\begin{equation}
\label{def-vert-mod-1} Y^M: V \to \Fields(M),\;  a \to Y^M(a,z) =
\sum_{n\in \ZZ} a^M_{(n)}z^{-n-1},
\end{equation}
 that satisfy the following axioms:
\begin{enumerate}
  \item $\vac^M (z)   = \id_M  $
  \item (Borcherds identity)
  \begin{align}\label{def-vert-mod-2}
     &\sum\limits_{j\geq 0} {m \choose j} (a^{}_{(n+j)} b
     )^M_{(m+k-j)}\\
  = &\sum\limits_{j\geq 0} (-1)^{j} {n \choose j}\{ a^M_{(m+n-j)} b^M_{(k+j)} - (-1)^{n}b^M_{(n+k-j)} a^M_{(m+j)} \}\nonumber
  \end{align}
\end{enumerate}

Note that we have unburdened the notation by letting
\[
a^M(z)=Y^M(a,z).
\]

A module $M$ over a graded vertex algebra $V$ is called {\em graded}
if $M = \oplus_{n\geq 0} M_n$ with
 $v_{k}M_l  \subset M_{l-k}$  (assuming $M_{n} = 0$ for negative $n$).

 A {\em morphism of modules} over a vertex algebra $V$ is a map $f: M \to N$
 that satisfies $f(v^M_{(n)}m) = v^N_{(n)}f(m)$ for $v\in V$, $m\in M$.
$f$ is {\em homogeneous} if $f(M_k)\subset N_k$ for all $k$.

\subsection{Examples}

\subsubsection{Affine vertex algebras.}
\label{Affine vertex algebras.} Let $\fg$ be a semisimple Lie
algebra and $\qf{}{}: S^2 \fg \to \pole$ an invariant  form on
$\fg$. The {\em affine Lie algebra} $\ghat$  associated with $\fg$
and $\qf{}{}$ is a central extension of $\fg\tensor \pole[t,t^{-1}]$
defined as follows. As a vector space, $\ghat = \fg\tensor
\pole[t,t^{-1}] \oplus \pole K$ and the Lie bracket is
$$ [x\tensor t^n, y\tensor t^m ] = [x,y] \tensor t^{m+n} + n \delta_{n+m, 0} \qf{x}{y}K$$
$K$ is a central element.

We denote $x\tensor t^n$ by $x_n$ and write $x(z) = \sum x_n
z^{-n-1}$.

 Let $\fg = \fn_{-} \oplus \fh \oplus \fn_{+}$ be a Cartan decomposition of $\fg$.

Denote $\ghat_{<} = \fg\tensor t\pole[t]$, $\ghat_{>} = \fg\tensor
t^{-1}\pole[t^{-1}]$ and $\ghat_{\leq} = \fg\tensor \pole[t] \oplus
\pole K$.

Define $\ghat_{+} = \fn_{+}\oplus \ghat_{>}$, $\ghat_{-} =
\fn_{-}\oplus \ghat_{<}$. Then $\ghat = \ghat_{+} \oplus \fh \oplus
\pole K \oplus \ghat_{-}$.

%
%{\em Weyl modules.}
% If  $\lambda$ is a weight of $\fg$ and $k\in \pole$ we define the Weyl module $V_{\lambda, k}$ to be the induced module $V_{\lambda, k} = \Ind^{\ghat}_{\ghat_{\leq}} \pole_{\lambda, k}$ where $\pole_{\lambda, k}$ is a 1-dimensional $\ghat_{\leq}$-module generated by a vector $v_{\lambda, k}$ such that
% $\ghat_{<}\vlk = 0$,           $h_0\vlk = \la(h) \vlk$
% for   $x\in \fg$ and $K\vlk =  \vlk$.

The space of invariant forms is one-dimensional, and we will
 let $\qf {}{}$ be that form  for which
 %the longest root has length $\sqrt{2}$.
 $(\theta, \theta) = 2$ where $\theta$ is the longest root.

 Introduce the following induced module
 \begin{equation}\label{def-of-aff-vert-alg}
 V_{k}(\fg)  = \Ind^{\ghat}_{\ghat_{\leq}} \pole_{k},
 \end{equation}
  where $\pole_{ k}$ is a 1-dimensional $\ghat_{\leq}$-module generated by a vector $v_{ k}$ such that
 $\ghat_{<}v_k = 0$,           $\fg v_k = 0$  and $Kv_k =  v_k$.

$V_{k}(\fg)$ carries a vertex algebra structure that is defined by
assigning to $x_{-1}v_k$ $x\in \fg$,  the field  $x(z) = \sum x_n
z^{-n-1}$. These fields generate $V_{k}(\fg) $.

 $V_{k}(\fg) $ is a graded vertex algebra with generators having conformal weight 1.
 For example,
 \begin{equation}
 \label{comp-o-1-aff}
 V_{k}(\fg)_0=\pole_k,\;V_{k}(\fg)_1=\fg\otimes t^{-1}v_k.
 \end{equation}

%%
%%
%%    If $W$ is a vector space with a symmetric bilinear form $(\cdot ,\cdot):S^2 W \to \Cplx$, define the Heisenberg Lie
%%

\subsubsection{Commutative vertex algebras.}
\label{Commutative_vertex_algebras} A vertex algebra is said to be
  {\em commutative} if $a_{(n)} b =0$ for $a$, $b$ in $V$ and $n\geq 0$.
   The structure of a commutative vertex algebras is equivalent to one %that
 of commutative associative algebra with a derivation.

If $W$ is a vector space we denote by $H_W$ the algebra of
differential polynomials on $W$. As an associative algebra it is a
polynomial algebra in variables $x_i$, $\de x_i$, $\de^{(2)}x_i$,
$\dots$ where $\set{x_i}$ is a basis of $W^*$. A commutative vertex
algebra structure on $H_W$ is uniquely determined by attaching the
field
 $x(z) = e^{z\de}x_i$  to $x\in W^*$.

 $H_W$ is equipped with grading such that
 \begin{equation}
 \label{weights-0-1-diff-poly}
 (H_W)_0=\pole,\; (H_W)_1=W^*.
 \end{equation}

\subsubsection{Beta--gamma system.}
\label{beta-gamma-system} Define the Heisenberg Lie algebra  to be
the algebra with generators $a^i_n$, $b^i_n$, $1\leq i\leq N$ and
$K$ that satisfy $[a^i_m, b^j_n] = \delta_{m, -n} \delta_{i,j} K$,
$[a^i_n, a^j_m] = 0$, $[b^i_n, b^j_m] = 0$.

Its Fock representation $M$ is defined to be the module induced from
the  one-dimensional representation $\Cplx_1$ of its subalgebra
spanned by $a^i_n$, $n\geq 0$, $b^i_m$, $m>0$ and $K$ with $K$
acting as identity and all the other generators acting as zero.

The {\em beta-gamma system} has $M$ as an underlying vector space,
the
 vertex algebra structure being determined  by
assigning the fields
$$a^i(z) = \sum a^i_n z^{-n-1}, \ \ b^i(z) = \sum b^i_n z^{-n}$$
to $a^i_{-1}1$ and $b^i_01$ resp., where $1\in \Cplx_1$.

This vertex algebra is given a grading so that the degree of
operators $a^i_n$ and $b^i_n$ is $n$. In particular,
\begin{equation}
\label{weight-0-1-b-g} M_0=\pole[b_0^1,...,b_0^N],\;
M_1=\bigoplus_{j=1}^{N}(b^j_{-1}M_0\oplus  a^j_{-1}M_0).
\end{equation}

\subsection{Vertex algebroids}
\label{Vertexalgebroids}

\subsubsection{Definition.}\label{def-of-vert-alg}Let  $V$ be a graded vertex algebra. We briefly recall
from \cite{GMS} basic results on the structure that is induced by
vertex operations on the subspace $V_{\leq 1} = V_0  + V_1$.

Let us define a 1-truncated vertex algebra to be a sextuple
 $(V_0\oplus V_1, \vac, \de, \opm, \zero, \one )$ where the operations $\opm, \zero, \one$ satisfy all the axioms
   of a vertex algebra that make sense upon restricting to the subspace $V_0 + V_1$. (The precise definition can be found in
    \cite{GMS}).
   The category  of 1-truncated vertex algebras will be denoted $\cV ert_{\leq
   1}$.

The notion of a 1-truncated vertex algebra is equivalent to that of
a vertex algebroid. For the  definition the reader is referred to
\cite{GMS}; in this note we only recall the main ingredients
and properties of a vertex algebroid.

For a graded vertex algebra $V$, set $A = V_0$, $\Omega = A\opm \de
A$ and $T=V_1 / \Omega$. The axioms of vertex algebra imply the
following:
\begin{enumerate}
\item
 $A = V_0$ is a commutative associative algebra with respect to ${}_{(-1)}$;
\item
 $\Omega$ is an $A$-module via $a\cdot \omega  =  a\opm \omega$ and  the translation map $\de:A \to \Omega$ is a derivation;
\item
 $T=V_1 / \Omega$ is a Lie algebra with bracket $\zero$ and a left  $A$-module via $\opm$;
  \item
$\Omega$ is a $T$-module with the action induced by $\zero$;  %%  $ T \times \Omega \to \Omega$.

 \item
the map $\zero: T \times A \to A$ defines an action of $T$ on $A$ by
derivations.
\item
the maps $_{(1)}: T\times\Omega\rightarrow A$ and $_{(1)}:
\Omega\times T \rightarrow A$ are $A$-bilinear pairings that satisfy
$\tau_{(1)}\omega=\omega_{(1)}\tau$ and are determined by
$\tau_{(1)}\partial a = \tau_{(0)}a$.
\end{enumerate}
The gadget (1)--(6) is quite classical; in particular, (1,5) mean
that $T$ is an $A$-{\it Lie algebroid} \cite{BB2}. Altogether (1--6)
were called an {\it extended Lie algebroid} in \cite{GMS}, a concept
that is equivalent to that of a {\it Courant algebroid}  --- this is
a remark of P.~Bressler, \cite{Bre}.

All of the vertex algebra structure on $V_{0}+ V_1$ comprises more
data than (1--6), but not much more. A {\em vertex algebroid} is a
quintuple
 $(A\oplus T\oplus \Omega, \de, \gamma, \langle, \rangle, c )$ where
$(A, \Omega, T, \de$ are as in (1--5),
  $\gamma: A \times T \to \Omega$ is a bilinear map,
$$\langle, \rangle : (T\oplus \Omega) \times (T\oplus \Omega) \to A  $$
is a symmetric bilinear pairing, and
$$
c: T\times T \to \Omega
$$
 is a
skew-symmetric bilinear
 pairing. These data satisfy a list of axioms to be found in \cite{GMS}.
 We will not record those axioms here -- they are a result of
 writing  down the restriction of the Borcherds identity to
 conformal weights 0 and 1 subspaces -- but we will supply the
 reader with a short dictionary:

Fix a splitting $V_1=T\oplus\Omega$, see item (3) above.

The map $\gamma$ is determined by the classical data (1--6), the
splitting chosen, and the
 Borcherds identity.

 The pairing $\langle,\rangle$ is an extension of the pairing from item
 (6); the extra part is
\begin{equation}
\label{defofpairing}
\langle\xi,\eta\rangle=\xi_{(1)}\eta,\;\xi,\eta\in T.
\end{equation}
 The map $c$, the key to extending the classical data (1--6) to a
 vertex algebroid is defined by (in the presence of the splitting)
 \begin{equation}
 \label{defofc}
\xi_{(0)}\eta=[\xi,\eta]+c(\xi,\eta),\;\xi,\eta\in T.
\end{equation}

The category of vertex algebroids, to be denoted $\cA lg$, is
defined in an obvious manner and is immediately seen to be
equivalent to that of 1-truncated vertex algebras.

\subsubsection{Truncation and vertex enveloping algebra functors}

There is an obvious  truncation functor $$t: \cV ert \to \cV
ert_{\leq 1} $$ that assigns to every vertex algebra a 1-truncated
vertex algebra.
  This functor admits a left adjoint \cite{GMS}
  $$ u:  \cV ert_{\leq 1}  \to \cV ert $$
   called a {\em vertex enveloping algebra functor}.

   In the context of vertex algebroids, these functors become  $\cA: \cV ert \to \cA lg$
  and its left adjoint $ U: \cA lg \to \cV ert .$

  \subsubsection{Examples.}
  \label{examples-of-vert-alg}
  Various examples of vertex algebras reviewed above are, in fact,
  vertex enveloping algebras of appropriate vertex algebroids:
  \begin{itemize}
  \item  in the situation of sect. \ref{Affine vertex algebras.},
  $\pole_k\oplus \fg\otimes t^{-1}v_k$ is a vertex algebroid, see
  (\ref{comp-o-1-aff}),
  and  $V_k(\fg)=U(\pole_k\oplus \fg\otimes t^{-1}v_k,...)$;

  \item in the situation of sect. \ref{Commutative_vertex_algebras},
  $\pole\oplus W^*$ is a vertex algebroid -- obviously commutative, see (\ref{weights-0-1-diff-poly}),
  and $H_W=U(\pole\oplus W^*)$;

  \item in the situation of sect.\ref{beta-gamma-system},
  $M_0\oplus M_1$ is a vertex algebroid, see (\ref{weight-0-1-b-g}),
  and $M=U(M_0\oplus M_1)$.
  \end{itemize}
  If we have two vertex algebras, $V$, $W$, then their tensor
  product $V\otimes W$ carries a vertex algebra structure defined as
  usual, see e.g. \cite{K,FBZ}, by letting
  \begin{equation}
  \label{def-of-tens-prod-algebra-vert}
  (v\otimes
  w)_{(n)}a\otimes b=\sum_{j=-\infty}^{+\infty}v_{(j)}a\otimes w_{(n-j-1)}b.
  \end{equation}
  One similarly defines the tensor product of two vertex algebroids.
  In the more convenient language of 1-truncated vertex algebras, if
   $V=V_0\oplus V_1$, $W=W_0\oplus W_1$ are two 1-truncated vertex
   algebras, then we define
   \begin{equation}
   \label{def-of-tens-prod-algebroid-vert}
   V\stackrel{\bullet}{\otimes}W=(V_0\otimes W_{0})\oplus( V_0\otimes
   W_{1}\oplus V_1\otimes W_{0}),
   \end{equation}
   \begin{equation}
   \label{def-of-oper-in-tens-prod-algebroid}
  (v\otimes
  w)_{(n)}a\otimes b=\sum_jv_{(j)}a\otimes w_{(n-j-1)}b,
  \end{equation}
  where unlike (\ref{def-of-tens-prod-algebra-vert}) the summation
  $\sum_j$
  is extended to those $j$ for which it makes sense.

  The vertex algebroid $M_0\oplus M_1$ will give rise to the
  simplest example of an algebra of chiral differential operators,
  the subject to which we now turn, and the tensor product
  $(M_0\oplus M_1)\stackrel{\bullet}{\otimes}
  (\pole\oplus W^{*})$ will be similarly used to construct an algebra of
  twisted chiral differential operators in sect.
  \ref{A_universal_twisted_CDO}; here $\pole\oplus W^{*}$ is a
  commutative vertex algebroid from sect.
  \ref{examples-of-vert-alg}.

\subsection{Chiral differential operators}
 \label{Chiral_differential_operators} A vertex
algebra $V$ is called an {\em algebra of chiral differential
operators } over $A$, CDO for short, if $V$ is the vertex envelope
of a vertex algebroid $\cA = A\oplus T\oplus \Omega$ such that $T =
\Der A$ and $\Omega  = \Omega^1_{A}$, the module of K\"ahler
differentials.

An algebra of chiral differential operators over $A$ does not exist for any $A$,
but it does exist locally on $Spec A$.

 To be more precise, a smooth affine variety $U= \Spec A $ will
  be called {\em suitable for chiralization}
 if $Der(A)$ is a free $A$-module admitting an abelian frame $\{\tau_1,...,\tau_n\}$.
 In this case
  there is a CDO over $A$, which is uniquely determined by the
  condition that
  $(\tau_i)_{(1)}(\tau_j)=(\tau_i)_{(0)}(\tau_j)=0$; in other words we let the ``quantum
  data'',
  $\langle,\rangle$ and $c$ vanish on the basis vector fields, cf. (\ref{defofpairing},\ref{defofc}).
Denote this CDO by $D^{ch}_{U,\tau}$.

\begin{thm}\label{class-cdo-local} Let $U=\Spec A$ be suitable for
chiralization with a fixed abelian frame $\{\tau_i\}\subset Der A$.

(i) For each closed 3-form $\alpha\in\Omega^{3,cl}_A$ there is a CDO
over $A$ that is uniquely determined by the conditions
\[
(\tau_i)_{(1)}\tau_{j}=0,\;(\tau_i)_{(0)}\tau_{j}=\iota_{\tau_i}\iota_{\tau_j}\alpha.
\]
Denote this CDO by $\cD_{U,\tau}(\alpha)$.

(ii) Each CDO over $A$ is isomorphic to $\cD_{U,\tau}(\alpha)$ for
some $\alpha$.

(iii) $\cD_{U,\tau}(\alpha_1)$ and $\cD_{U,\tau}(\alpha_2)$ are
isomorphic if and only if there is $\beta\in\Omega^{2}_A$ such that
$d\beta=\alpha_1-\alpha_2$. In this case the isomorphism is
determined by the assignment
$\tau_i\mapsto\tau_i+\iota_{\tau_i}\beta$.
\end{thm}

If $A=\pole[x_1,...,x_n]$, one can choose $\partial/\partial x_j$,
$j=1,...,n$, for an abelian frame and check that the beta-gamma
system $M$ of sect. \ref{beta-gamma-system} is a unique up to
isomorphism CDO over $\pole^n$. A passage from $M$ to Theorem
\ref{class-cdo-local} is accomplished by the identifications
$b^j_01=x_j$, $a^j_{-1}1=\partial/\partial x_j$.

The construction of CDOs from Theorem \ref{class-cdo-local} can be
sheafified; however, what one gets is not a sheaf but rather a gerbe
over a smooth variety $X$ bound by the complex
$\Omega^2_X\rightarrow\Omega^{3,cl}_X$. The existence of global
objects in this gerbe depends on vanishing of a certain
characteristic class of $X$. For the precise description of the
situation we refer the reader to \cite{GMS}, here let us just note
that for a smooth variety $X$ with vanishing first Pontrjagin class
there exist such sheaves; they are called {\em sheaves of chiral
differential operators.}

Let $\cdo$ denote any of such sheaves.   This is a graded sheaf. A
straightforward consequence of the construction is that
 $(\cdo)_0 \iso \oX$, $\Omega^1_X \subset (\cdo)_1$ and  $(\cdo)_1 / \Omega^1_X \iso \cT_X$.

\bigskip

Another "classical"  object that we  attach to a vertex algebra $V$
is the universal enveloping algebra $U_A T$ of
 the $A$-Lie algebroid $T$.
By definition,  $U_A T$ is the quotient of the tensor algebra
 \begin{center} $Tens(A\oplus T)= \oplus_{i\geq 0} (A\oplus T)^{\otimes i}$
 \end{center}
  modulo the ideal $R$ generated by the  elements $a\otimes b -  ab$,
  $\tau \otimes a- a \otimes \tau - \tau (a)$, $\tau\otimes\xi-\xi\otimes\tau - [\tau,\xi]$, $a\otimes \tau -  a\tau$, $1_A  - 1_\Cplx$.
In the next section we will see how this algebra appears via Zhu's
construction.

\section{Zhu's correspondence}

The work  \cite{Zhu} revealed a beautiful and nontrivial connection
between the world of vertex algebras and that of associative
algebras. The main result of Zhu's theory states that to each graded
vertex algebra $V$ one can naturally attach an associative algebra,
to be denoted $Zhu(V)$, such that there is a one-to-one
correspondence between simple $Zhu(V)$-modules and simple
$V$-modules.

The aim of this section is to prove the theorem below. This is a
statement connecting Zhu algebra of $V$ to the universal enveloping
algebra $U_A T$ defined in the previous section. The basic
observation is that there is a natural associative algebra morphism
 \begin{equation}
  \label{map-alpha}
 \alpha: U_A T \to  Zhu(V),
 \end{equation}
  to be constructed in Subsection \ref{subsec_alpha}.
\begin{thm} \label{zhu-envelope}
(1)  If $V$ is generated by $V_0 + V_1$, then the map $\alpha$ is
surjective.

(2) If $V$ is a vertex enveloping algebra of $V_0+V_1$ and $T$ and
$\Omega$ are  free $A$-modules, then $\alpha$ is an isomorphism.
\end{thm}

\begin{rem}
  If $V = V_{k}(\fg)$, see sect.\ \ref{Affine vertex algebras.},
  then $U_A T = U\fg$, and Theorem \ref{zhu-envelope} follows from the
  isomorphism
  $$
  U\fg \iso Zhu(V_k (\fg))
  $$
  established in \cite{FZ}.
  It is fair to say that Theorem \ref{zhu-envelope} is a variation on the theme of \cite{FZ}.
\end{rem}

As a corollary, we will have a description of the sheaf of Zhu
algebras for the vertex algebra of twisted differential operators.

The proofs of many of the auxiliary results  below can be found
elsewhere, e.g. \cite{Zhu,FZ,R,MZ,L}.

\subsection{Definition of the  Zhu algebra}
\subsubsection{Motivation}     %\marginpar{subsection slightly rewritten}
\label{Motivation-zhu} While vertex operations ${}_{(n)}$, $n\in
\ZZ$ satisfy axioms as remote from associativity as Borcherds
identity, the
 endomorphisms $v_{(n)}$ belong to an associative algebra $\End V$.
 Furthermore, if
 $M$ a graded $V$-module,  then there are maps
 $$
 V\rightarrow \End(M),\; v\mapsto v^M_0
 $$
 and by restriction
$$
 V\rightarrow \End(M_0),\; v\mapsto v^M_0|_{M_0}
 $$

 Now one can ask if there is an operation $*$ on $V$ that makes the
 latter an  algebra morphism for any $M$. The answer is yes, and
 in order to find such an operation let us look at the Borcherds identity for a $V$-module $M$:
$$   \sum\limits_{j\geq 0} {m \choose j} (a\ops{n+j} b )^M_{(m+k-j)}  = \sum\limits_{j\geq 0} (-1)^{j} {n \choose j}\{ a^M_{(m+n-j)} b^M_{(k+j)} - (-1)^{n}b^M_{(n+k-j)} a^M_{(m+j)} \}  $$
in the conformal weight notation
$$
   \sum\limits_{j\geq 0}
    {m \choose j} (a_{(n+j)} b )^M_{m+k-\Delta_a - \Delta_b +n+2}
   =
    \sum\limits_{j\geq 0}
     (-1)^{j} {n \choose j}\{ a^M_{m+n-j -\Delta_a +1} b^M_{k+j -\Delta_b +1}
     $$
     $$
    - (-1)^{n} b^M_{n+k - j -\Delta_b +1} a^M_{m+j -\Delta_a +1} \}
$$
and consider the case when $m = \Delta_a$, $n = -1$,  and both sides
are degree 0 morphisms, which requires $k = \Delta_b -1$. We obtain
$$   \sum\limits_{j\geq 0} {\Delta_a \choose j} (a_{(-1+j)} b )^M_{0} =
\sum\limits_{j\geq 0} ( a^M_{-j } b^M_{j }  +  b^M_{- j - 1} a^M_{j  +1} )  $$
Restricting this to $M_0$ will give us
$$   (\sum\limits_{j\geq 0} {\Delta_a \choose j} a_{(-1+j)} b   )^M_{0}|_{M_0} = a^M_0 b^M_0|_{M_0} $$
which means that for the desired  operation we can take the
following
$$a*b = \sum\limits_{j=0}^{\Delta_a} { \Delta_a \choose
                                            j} a_{(-1+j)}b.$$

This operation is not associative. However, it is shown in Zhu's
work  that there is a subspace $O(V) \subset V$ that is an ideal
with respect to this operation and acts as zero
 on $M_0$ for each $M$,  and such that the induced multiplication on $V/O(V)$ is associative.
Specifically, $O(V) = (\de  + H)V * V $, where $Hv = \Delta_v v$ for
homogeneous $v$.
 It is straightforward to verify that $v_0 = 0$ for $v\in O(V)$, as $(\de v)_0 = - \Delta_v v_0$.
 What is more remarkable is that $O(V)$ is an ideal with respect to $*$ and that $*$ is associative modulo
 $O(V)$. Furthermore, the associative algebra $(V/O(V),*)$ carries
 some
 essential information on the category of $V$-modules.

\subsubsection{Formal definition.}
\begin{defn} For homogeneous $a\in V$ define the Zhu multiplication

$$a*b = \sum\limits_{i=0}^{\Delta_a} { \Delta_a \choose
                                            i} a_{(i-1)}b.$$

More generally, for  $n\in\ZZ$ define  %%other Zhu operations:

$$a*_n b = \sum\limits_{i=0}^{\Delta_a} {
                                            \Delta_a \choose
                                            i
                                          } a_{(n+i)}b $$
\end{defn}

 To make this operation associative one has
to pass over to a properly chosen quotient of $(V, *)$.

 Denote  $O(V)=V*_{-2} V$. One can show that $O(V) = V*dV = dV
* V$ where $d= \de + H$\

One has the following
\begin{prop} (1) $V*O(V)\subset O(V)$, $O(V)*V\subset O(V)$

(2) $a*(b*c)- (a*b)*c \in O(V)$ for all $a$, $b$, $c$ in $V$.
\end{prop}

\begin{defn} Define the Zhu algebra to be the space $$Zhu(V)=V/O(V)$$ endowed with the
multiplication induced by $*$.
\end{defn}

It follows from the proposition above that the  Zhu algebra is an
associative algebra. It is naturally a filtered algebra with the
filtration induced by the conformal weight filtration of the vertex
algebra $V$. Specifically, $F^m Zhu(V) =
\pi(\bigoplus\limits_{i=0}^{m} V_i)$ where $\pi: V \to V/ O(V)$ is
the natural projection map.

Let $\cV ert$ denote the category whose objects are $\ZZ_{\leq
0}$-graded vertex algebras and the morphisms are graded vertex
algebra maps.

The  correspondence $V\mapsto Zhu(V)$ provides a functor from $\cV
ert$ to the category of filtered associative algebras.

\subsection{Relation of $Zhu(V)$ to $U_A T$.}

\subsubsection{}\label{subsec_alpha}

Recall that $A=V_0$ is an associative commutative algebra with
multiplication $\opm$ and $T =V_1/ V_0 \opm \de V_0$ is an $A$-Lie
algebroid.

Since $\Omega = A\opm \de A = A * dA$ is a subset of  $V*dV = O(V)$,
we have a natural map $\bar{\alpha}: T = V_1 / \Omega \to V/ O(V) =
Zhu(V)$. We denote $\bar{\alpha}(\tau)=\btau$.

\begin{lem} $A$ is naturally embedded into $Zhu(V)$.
\end{lem}
\noindent \textbf{Proof.} First, notice that elements of the form
$dw*v$ do not have  a degree 0 component. Indeed, the lowest
conformal weight summand in $dw*v \in \oplus_{i =
\Delta_v}^{\Delta_v + \Delta_w +1 } V_i$ is equal to  $(dw)_0 v = 0$
(since $d = \de + H$ and  $(\de w)_0   =  - \Delta_w w_0$).

 Thus $V_0\cap O(V) = 0$
 and the restriction of the projection $\pi :V \to Zhu(V)$ to $V_0$  is injective.
  Since $a*b = a\opm b$ for $a,b \in V_0$, this is an algebra embedding. $\qed$

\begin{lem}
\label{alpha-defining-thm} The natural map $\bar{\alpha}: T \to
Zhu(V)$ extends to an algebra homomorphism $\alpha : U_A T \to
Zhu(V)$
\end{lem}
\noindent \textbf{Proof}. The inclusion $A\hookrightarrow Zhu(V)$
and the map $\bar{\alpha}: T\to Zhu(V)$ uniquely determine algebra
morphism $Tens(A\oplus T) \to Zhu(V)$.

We have an exact sequence
$$
0 \to R \to Tens(A\oplus T) \to U_A T \to 0,
$$
cf. the end of sect.\ \ref{Chiral_differential_operators}

 Under the map $Tens(A\oplus T) \to Zhu(V)$ the generators of $R$ are mapped to
$R_1=a*b- a\opm b$,
 $R_2= \btau*a- a*\btau - \tau (a)$,
   $R_3=\btau*\bxi-\bxi*\btau - \overline{[\tau,\xi]}$,
 $ R_4= a*\btau  - a\opm\btau$, and $1-1$.
To finish the proof, it suffices to show that $R_i=0$, $i=1, 2, 3,
4$.

1)  $R_1=0$ due to the algebra inclusion $A\hookrightarrow
 Zhu(V)$.

2). Let $t$ denote any lifting of $\tau$ to $V_1$. We  check that
$t*a - a*t - \tau(a) \in O(V)$

Recall that for $x\in V_1$ the Zhu operation reduces to $x*v = x\opm
v + x\zero v$ and by definition, $\tau(a) = t\zero a$. Hence
$$t*a - a*t - t(a) = t\opm a + t\zero a - a\opm t - t\zero a =$$
$$= t\opm a - a\opm t = [t\opm, a\opm]\vac = \de (t\zero a) \in \de(V_0)\subset  O(V)$$

3) Let $q$ be any lifting of $[\tau, \xi]$ to $V_1$. Then $q-t\zero
x \in \Omega$ for any lifting $x$ of $\xi$ and $t$ of $\tau$. Thus
it suffices to check
$$
t*x- x*t -t\zero x = t\opm x + t\zero x -x \opm t - x\zero t -
t\zero x =
$$
$$
 =[t\opm, x\opm]\vac -x\zero t = \de(t\zero x) - x\zero t =\de(t\zero x) + t\zero x - \de(t\one x) \in O(V)
 $$
4) $R_4=0$ follows from $a*v = a\opm v$ for any $v\in V$, $a\in A$.

Thus, the map factors through $U_A T$. $\qed$

\subsection{Properties of the Zhu algebra when $V$ is generated by
$V_{\leq 1}$}
\subsubsection{}
 Notice that
\begin{equation}
\label{a_star_n} a*_n v = a_{(n)}v,\ \ \ \ {\rm{for }}\ a\in V_0
\end{equation}
\begin{equation}
\label{x_star_n} x*_n v = x_{(n)}v + x_{(n+1)}v \ \ \ {\rm{for }}\
x\in V_1
\end{equation}

We will start with deriving a different presentation of the ideal
$O(V)$. First, we observe the following:

\begin{prop} (1) $V*_nV \subset V*_{n+1}V$ for all $n\neq -1$
\end{prop}

\textbf{Proof.} Straightforward, see e.g. [R] $\qed$

\bigskip
Consequently, $V*_{n} V \subset O(V)$ for all $n\leq -2$. In
particular, the subspace
\begin{equation}
%\nonumber
O'(V)= \spn\{ a_{(n)}v, (x_{(n)}+ x_{(n+1)})v | \ n\leq -2,\  a\in
V_0, x\in V_1, v\in V\}
\end{equation}
 lies in $O(V)$.

The aim of the next lemmas is to show that if $V$ is generated by
$V_0 +V_1$ then $O'(V)$ is in fact all of $O(V)$.

First, we notice that $O'(V)$ is invariant under the action of
$f_{m}$ for  $f\in V_{\leq 1}$, $m\leq 0$.

\begin{lem}\label{closure}   %  (Closure property of $O'(V)$).
We have
 $$b_{(m)} O'(V)\subset
O'(V), \ \  y_{(n)} O'(V)\subset O'(V)$$ for any  $b\in V_0$, $y\in
V_1$, $m\leq -1$, $n\leq 0 $.
\end{lem}

\textbf{Proof}. The space $O'(V)$ is spanned by  elements of the
form  $a\ops{k}v$ and $(x\ops{k}+x\ops{k+1})v$, $k\leq -2$. We need
to show that for $a,b\in V_0$ and $x,y \in V_1$ the elements $ A=
b_{(m)}a_{(k)}v, \ B=b_{(m)}(x_{(k)}+ x_{(k+1)})v, \
C=y_{(n)}a_{(k)}v, \ D=y_{(n)}(x_{(k)}+ x_{(k+1)})v $
 can also be written in such form.

Cases A and B. For $m\leq -2$ we have $b_{(m)}O'(V) \subset
b_{(m)}V\subset O'(V)$ by definition of $O'(V)$. For $m=-1$ we use
the commutator identities. In one case we have
$$
A=b\ops{-1}a_{(k)}v= a_{(k)} b\opm v + [b\opm, a_{(k)}]v = a_{(k)}
b\opm v \in O'(V),
$$
 in the other we use
$
 [b_{(-1)}, x_{(k+1)}]v =
(b_{(0)}x)_{(k)} v \in O'(V), \ k\leq -2 $
 to conclude that
$$
B=b\opm (x_{(k)}+ x_{(k+1)})v =  (x_{(k)}+ x_{(k+1)}) b\opm v\ +
(b\zero x)_{(k-1)} v + (b\zero x)_{(k)} v
$$
 is also in $O'(V)$.

The remaining cases are also implied by the commutator formula:
$$
C=y\ops{n} a\ops{k} v = a\ops{k} y\ops{n} v  + [y_{(n)}, a_{(k)}]v
=  a\ops{k} y\ops{n} v  + (y_{(0)}a)_{(n+k)}v  \in %% V_0 \ops{n}V \subset
O'(V)
$$
for $k \leq -2$, $n\leq 0$.
 $$
 D=
y\ops{n} (x_{(k)}+ x_{(k+1)})v= (x_{(k)}+ x_{(k+1)}) y\ops{n}v +
[y_{(n)}, x_{(k)} + x_{(k+1)}]v
$$
To show $D\in O'(V)$ we notice
$$
[y_{(n)}, x_{(k)} + x_{(k+1)}]v = ((y_{(0)}x)_{(n+k)}
+(y_{(0)}x)_{(n+k+1)})v +
$$
$$
+ n(y_{(1)}x)_{(n+k-1)}v + n(y_{(1)}x)_{(n+k)}v =
$$
$$
= (y_{(0)}x)*_{n+k}v + n(y_{(1)}x)*_{n+k-1}v + n(y_{(1)}x)*_{n+k}v,
$$
 all three terms in $O'(V)$.
$\qed$
\bigskip

From now on, $V$ is a vertex algebra generated by  $V_{\leq 1}$,
unless stated otherwise.

 \begin{lem}  $O(V) = O'(V)$.
\end{lem}

\textbf{Proof} We need to show that $u*_{-2}v \in O'(V)$ for all
$u,v\in V$. Following the proof of Rosellen (\cite{R}, Proposition
6.2.5) we show that $u*_{n}v \in O'(V)$ for all $n\leq -2$ by
induction on $\Delta_u$.

Basis of induction: $\Delta_u =0$, that is,  $u\in A$. Then
$u*_{n}v\in O'(V)$ for all $n\leq -2$ by definition of $O'(V)$.

% Induction step.
Suppose we showed that $u*_{n}v \in O'(V)$ for all $u, v\in V$ such
that $\Delta_u =k$ and all $n\leq -2$. Now we need to prove it for $V_{k+1}$. Any
element of $V_{k+1}$ is a sum of elements of the form $x_{(-r)}u$ or
$b_{(-r-1)}u$ where $\Delta_{u} \leq k$, $r\geq 1$.

Consider  $b_{(-r-1)}u$ with $b\in V_0$. Since $b_{(-r-1)}=
b*_{-r-1}$, we can use the associativity formula (see [R,
Proposition 6.2.2])
 $$(b*_{-r-1}u)*_{-2-t}v  =
\sum\limits_{i,j=0}^{\infty} (-1)^{i}{r \choose j}{-r-1 \choose i} (
b*_{-r-1-i}(u*_{-2-t+i+j}v) $$
$$- (-1)^{-r-1} u*_{-3-t-r-i+j}({b*_{i} v}) ) )$$
The term $b*_{-r-1-i}(u*_{-2-t+i+j}v)$ is in $O'(V)$ since
$-r-1-i\leq -2$.

For the second term, we can assume that $0\leq j\leq r$ since
otherwise ${r \choose j}=0$ (as $r\geq 1$). Then $-3-t-r-i+j\leq
-3-t-i\leq -3$ and the term is in $O'(V)$ by the induction
assumption.

The proof for the case $ x_{(m)}u$ repeats Rosellen's proof in
Proposition 6.2.5.
%%%%      \marginpar{ Citing the book...}
For the sake of completeness, we reproduce it here, notation
slightly changed.

We need to show that $(x_{(m)}u)*_n v \in O'(V)$ for all $m \leq
-1$, $n\leq -2$. By (\ref{x_star_n}) we have
\begin{equation}
  (x_m u )*_n v  = (x*_m u)*_n v - (x_{m+1} u)*_n v.
\end{equation}
The last term of r.h.s. is in $O'(V)$ by the induction assumption.

By [R], Proposition 6.2.2 we have
$$
(x*_m u)*_n v = \sum\limits_{i,j\geq 0} (-1)^{i} {-m-1 \choose j}{m
\choose i} (x*_{m-i} (u*_{n+i+j} v) - (-1)^m  u*_{n+j+m-i} (x*_i
v))
$$
By induction, $
  u*_{n+j+m-i} (x*_i  v) \in O'(V)
$ since $j   \leq    -m-1$.

From (\ref{x_star_n}) it follows that $ x*_{m-i}(u*_{n+i+j}v) \in
O'(V) $ if $m-i \leq -2$,
 that is,
if $m\leq -2$ or $i>0$. If $m = -1$ and $i = 0$ then $j = 0$.
 By the induction assumption,
$u*_{n}v \in O'(V)$. Thus, using Lemma \ref{closure} we have
$x*_{-1} (u*_n v )  = x\opm  (u*_n v ) + x\zero  (u*_n v )  \in
O'(V)$. The lemma is proved. $\qed$

% The next lemma gives the ideal $O(V)$ a more suitable form for the proof of  Thm
% 2.

\begin{cor}
\label{cor-O-presentation} $O(V) = \spn \{ \omega_{(n+1)} v,\
(x_{(n)} + x_{(n+1)})v , \  v\in V,\ n\leq -2,\ x\in V_1, \   \omega
\in \Omega \}$
\end{cor}

\textbf{Proof.}  Denote $$O''(V) = \spn \{ \omega_{(n+1)} v,\
(x_{(n)} + x_{(n+1)})v , \  v\in V,\ n\leq -2 \}.$$ We show that
$O'(V) = O''(V)$.

Clearly, $O'(V) \subset O''(V)$, since $a_{(n)}= -\frac{1}{n+1} (\de
a)_{(n+1)}$, $n\leq -2$.

Now check $O''(V) \subset O'(V)$.  Let $\omega = a \de b$, $a,b \in
A$ and let $n\leq -1$. Then $$\omega\ops{n}v = \sum\limits_{j\geq 0}
(\de b)\ops{n-1-j}a\ops{j}v  + a\opm (\de b)\ops{n}v +
\sum\limits_{j\geq 1} a\ops{-1-j} (\de b)\ops{n+j}v =$$

$$ = \sum\limits_{j\geq 0}
 (j+1-n)b\ops{n-2-j}a\ops{j}v   - n a\opm b\ops{n-1}v +
\sum\limits_{j\geq 0} a\ops{-2-j} (\de b)\ops{n+1+j}v =$$

Both sums clearly are in $O'(V)$. The middle term is in $O'(V)$
since $b\ops{n-1}v\in O'(V)$ and $a\opm O'(V)\subset O'(V)$, see
Lemma \ref{closure}. $\qed$

\subsubsection{}
Fix a vertex algebra $V$  which is generated by $V_{\leq 1}$.

Using results obtained above we can state the following two lemmas.

\begin{lem}\label{zhuspan}  $Zhu(V)$ is spanned by $1$ and
$$  \pi(a^{}_{(-1)}x^1_{(-1)}\dots x^k_{(-1)}\vac),\
\textrm{where}\ \ a\in V_0,\ x^i\in V_1, \ 1\leq i\leq k, \ k\geq
0$$ where $\pi: V\to V/O(V)$ denotes the projection map.
\end{lem}

\textbf{Proof}. Indeed, if $V$ is generated by $V_{\leq 1}$, then
$V$ is spanned by $V_0$ and  monomials of the form $a^{}_{(m)}
x^1_{(-p_1-1)}\dots x^k_{(-p_k -1)}$, $a\in V_0$, $x^i \in V_1$. All
such monomials with $m\leq -2$ are by definition in $O'(V)$, so they
have no contribution to $Zhu(V)$.

Since  $(x\ops{n-1} + x\ops{n})v \in O'(V)$, $n\leq -1$ and
$x\ops{n} O'(V) \subset O'(V)$ we can conclude that
$$a^{}_{(-1)} x^1_{(-p_1-1)}\dots x^k_{(-p_k -1)}\vac \equiv
(-1)^{p_1} a^{}_{(-1)} x^1_{(-1)} x^2_{(-p_2-1)}\dots x^k_{(
-p_k-1)}\vac \equiv $$
$$ \equiv (-1)^{p_1 + p_2} a^{}_{(-1)} x^1_{(-1)} x^2_{(-1)}\dots x^k_{(
-p_k-1)}\vac \equiv \dots \equiv (-1)^{p_1 +\dots p_k} a^{}_{(-1)}
x^1_{(-1)}\dots x^k_{(-1)}$$ modulo $O'(V)$ $\qed$

\begin{lem}{\label{zhuspan2} } $Zhu(V)$ is generated
by the image of $V_{\leq 1}$.
\end{lem}

\textbf{Proof}. In view of the previous lemma, we just need to show
that the elements $\pi(a^{}_{(-1)}x^1_{(-1)}\dots x^k_{(-1)}\vac)$
are products of elements of $\pi(V_{\leq 1})$.
%%This is true indeed, due to another lemma in Zh's article
%
%Let us show that
% $$a^{}_{(-1)}x^1_{(-1)} \dots x^r_{(-1)}\vac  \equiv a*x^r *\dots x^1 \mod
% O(V).$$
We have  $a*v = a_{(-1)}v$ from definition and $v*x \equiv x_{(-1)}v
\mod O(V)$ from \cite{Zhu},  Lemma 2.1.3.  Therefore
$$x^r * \dots *x^2 * x^1 \equiv  x^1 \opm (x^r * \dots * x^2) \equiv \dots \equiv x^1 \opm x^2 \opm \dots x^r $$
and thus $ a*x^r * \dots *x^2 * x^1 \equiv  a\opm x^1 \opm x^2 \opm
\dots x^r \mod O(V)$. $\qed$

\bigskip

\subsubsection{Filtration of the Zhu algebra.} Here we briefly recall
different filtrations of $V$ that were dealt with in  \cite{GMS} and
the corresponding induced filtrations on $Zhu(V)$.

 Let $V$ be a $\ZZ_{\geq 0}$-graded vertex algebra generated by its first two components.
There is an obvious {\em  conformal weight filtration} $\cF = \{ F^n
V , n\geq 0\}$ defined by
$$F^n V = \bigoplus\limits_{i=0}^{n} V_i$$

In addition, there is a natural  filtration
  $\cG  = \{ G^n V, \ n \geq 0\}$
  by "number of vector fields".
By a vector field we mean an element of $T = V_1/\Omega$; by abuse
of language we will  call a vector field any element of $V_1$ that
projects onto a nontrivial element of $T = V_1/\Omega$.

This filtration is defined as follows:  the space $G^m V$ is the
space spanned by monomials $s\ops{n_1}^1\dots s\ops{n_r}^r \vac$,
where $s^i$ are elements either of $ V_0$ or of  $V_1$,  $n_i\leq
-1$, and  the number of vector fields among $s^1 , \dots , s^r$ is
less than or equal to $m$, i.e. $| \{  i : \ s_i\in V_1\backslash
\Omega    \}  |  \leq m$.

Clearly, $\cG $
 is an increasing exhaustive filtration of $V$.

\begin{lem} Filtration $\cG$ has the following property:
$$ G^{i}V \ops{n} G^j  V\subset G^{i+j}V\ \ \ \  {\textrm{for }}\ n\leq -1  $$
$$ G^{i}V \ops{n} G^j V \subset G^{i+j-1}V\ \ \ \  {\textrm{for }}\ n\geq 0   $$
\end{lem}

The proof is left as an exercise. An interested reader may find the proof of this fact in a more general setting in \cite{R} .

\begin{lem} Both $\cF$ and $\cG$ induce the same filtration on
$Zhu(V)$.
\end{lem}

\textbf{Proof}. We need to show $F^i V\subset G^i V + O(V)$ and $G^i
V \subset F^i V + O(V)$.

We have $F^i V\subset G^i V$ since a monomial of conformal weight
less than or equal to $i$ has at most $i$ vector field operators in
its formula.

Let us show $G^i V\subset F^i V + O(V)$. Any element of $V$ is a sum
of monomials of the form $a\opm x^1_{(n_1)} \dots x^s_{(n_s)} \vac$
with $a\in V_0$, $\{x^i\}_{ 1\leq i \leq s} \subset V_{1}$.

Let  $v\in G^i V$ be such a monomial. If $v$ contains  $x\ops{n}$
with $x\in\Omega$, then it is in $O'(V)$. Otherwise $s=i$ and, using
the proof of  Lemma \ref{zhuspan} one can show that  $v$ is equal to
$a\opm x^1_{(-1)} \dots x^s_{(-1)} \vac$ modulo $O'(V)$. This
monomial has conformal weight $i$, so $v\in F^i V + O(V)$.
 $\qed$

\bigskip

The enveloping algebra $U_A T$ is naturally a  filtered algebra by
(images of) $T^{\otimes n}$, $n\geq 0$. We have

\begin{lem} The map $\alpha: U_AT\rightarrow Zhu(V)$ of Lemma \ref{alpha-defining-thm} is a morphism of
filtered algebras.
\end{lem}

\textbf{Proof.} Let $F^i=F^i U_A T$ be the $i$-th filtration
subspace of $U_A T$. Clearly, $F^0= A$, $F^1=T $ and $F^k\subset
(F^1)^k$. Since $\alpha$ is a homomorphism and $\alpha(F^1) =
\alpha(T)\subset G^1 Zhu(V) $ we have $\alpha(F^k) \subset (G^1
Zhu(V))^k \subset G^k Zhu(V)$ $\qed$

\subsubsection{The Zhu algebra of enveloping algebras.
Proof of Theorem \ref{zhu-envelope}} Let  $V$ be a vertex algebra.
Recall that associated to $V$  is a vertex algebroid $\mathcal{A}_V
= (A, T, \Omega, \dots)$. In the subsection \ref{subsec_alpha} we
defined the map of associative algebras  $\alpha: U_A T \to  Zhu
(V)$. Theorem  \ref{zhu-envelope} states that $\alpha$ is surjective
when $V$ is generated by $V_0 + V_1 $ and it is an isomorphism when
$V$ is a vertex envelope of $V_0 + V_1$ and $\Omega$ and $T$ are
free $A$-modules. Now we are ready to complete the proof of Theorem
\ref{zhu-envelope}.

%  \textbf{Proof of Theorem \ref{zhu-envelope}}

(1) Surjectivity follows from Lemma \ref{zhuspan2}.

(2) We will show that the induced map
$$\bar{\alpha}:  Sym_A(T) \to \gr Zhu(V)$$
 is an isomorphism of commutative algebras.
%%   by factoring it into a product of two isomorphisms.
%%  $$Sym_A (T) \stackrel{\beta}{\longrightarrow} (\grH V) / \Symb O(V)
%%  \stackrel{\Phi}{\longrightarrow} \gr Zhu(V) $$

Recall \cite{GMS} that there are canonical filtrations $\cH$ on $V$
and $\cJ$ on $\grG V$ such that $\grH V = {\rm{gr}}^J \grG V $ and a
canonical map
\begin{equation}
\label{graded_iso_sym}
 \beta: Sym_A(\bigoplus\limits_{i\geq 0}T^{(i)} \oplus
\bigoplus\limits_{i\geq 0 }\Omega^{(i)})  \to  \grH V
\end{equation}
which is an isomorphism of commutative algebras provided that
$\Omega$ and $T$ are free $A$-modules.

From the definition of the filtration $\cH$ and Corollary
\ref{cor-O-presentation}, it follows that $I = \rm{Symb}^\cH O(V)$
is the ideal in $\grH V$ generated by symbols of $\de^{(i)}\omega$,
$i\geq 0$ and $\de^{(j)}\tau$, $i\geq 1$ where $\tau \in V_1$.

Hence we have a map
\begin{equation}
  \bar{\beta}:  Sym_A T  \iso Sym_A(\bigoplus\limits_{i\geq 0}T^{(i)} \oplus
\bigoplus\limits_{i\geq 0 }\Omega^{(i)}) / \beta^{-1}(I) \to \grH V
/ I
\end{equation}
which is a commutative algebra isomorphism.

It remains to notice that $\bar{\beta}$ is a composition of
$\bar{\alpha}$
with the natural  map    %%  isomorphism
$\gr Zhu(V) \iso \grG V/ \Symb^\calG O(V)  \to \grH V / I$ which
implies injectivity of $\alpha$. $\qed$

\begin{rem}
  The condition that $\Omega$ is a free $A$-module can be dropped, with a slight change to the proof.
\end{rem}

\bigskip

\subsection{ The Zhu algebra of a CDO}
 In this subsection we apply the results  obtained above  to the sheaf  of chiral differential operators.

If $\cV$ is a sheaf of vertex algebras on a variety $X$, we denote
$Zhu(\cV)$ the sheaf associated to the presheaf $U \mapsto
Zhu(\cV(U))$.

The Theorem \ref{zhu-envelope} has the following corollary:
\begin{cor} \label{zhuofcdo}
  Let $\cD^{ch}_X$ be a CDO on $X$. Then
 $$
 Zhu(\cD^{ch}_X) \iso \cD_X.
 $$
\end{cor}

For $X = {\mathbb{A}}^n$ this fact was proved in [L].

\noindent \textbf{Proof.} First, let us show that
 for  $U$ that is suitable for chiralization, sect.\ \ref{Chiral_differential_operators}, we have an isomorphism
 $Zhu(\cD^{ch}_U) \iso \cD_U$.

The algebra $\cdo(U)$ is an enveloping algebra of its 1-truncated
part, therefore by Theorem \ref{zhu-envelope} there is a natural
isomorphism $\alpha_U: \cD_X(U) = U_{\oX(U)}\calT(U) \to
Zhu(\cdo(U))$. For $V\subset U$ the isomorphisms $\alpha_V$ are
compatible so that we have an isomorphism of sheaves.

To show that $Zhu(\cdo)$ is globally isomorphic to $\cD_X$ it is
enough to notice that we have the embeddings $\oX=(\cdo)_0
\hookrightarrow Zhu(\cdo)$ and $(\cdo)_1 /\Omega^1 \hookrightarrow
Zhu(\cdo)$ that implies $\cT_X \hookrightarrow Zhu(\cdo)$. This is
due to the fact that transition functions for $\cdo$ are constructed
in \cite{GMS} in such a way that $(\cdo)_1 /\Omega^1$ is equal to
$\cT_X$. $\qed$

\subsection{The Zhu correspondence for modules.}

\begin{thm}(\cite{Zhu})
 Let  $M$ be a graded module over $V$. Then the top component $M_0$ is a module over the Zhu algebra $Zhu(V)$.
 The assignment $M\mapsto M_0$ establishes a 1-1 correspondence
 betweem isomorphism classes of irreducible graded $V$-modules
 and irreducible $Zhu(V)$-modules.
\end{thm}

Theorem \ref{zhu-envelope} implies the following result, which was
originally proved by other methods in \cite{LiYam}.

\begin{cor}
\label{graded-irrep-cdo} With the assumptions of Theorem
\ref{zhu-envelope} (2), there is a 1-1 correspondence betweem
isomorphism classes of irreducible graded $V$-modules
 and irreducible $U_AT$-modules.
 \end{cor}

\begin{rem}
It is enough to assume that $M$ is a filtered $V$-module.
\end{rem}
\begin{rem}
Let $V\rm{-}Mod$ denote the category of graded $V$-modules $M$. The
assignment
\begin{equation}
\label{top-component} \Phi: M \mapsto M_0
\end{equation}
is  a functor from the category of graded (resp.  filtered)
$V$-modules to the category of $Zhu(V)$-modules, with the obvious
action on maps.
\end{rem}

\subsubsection{Left  adjoint to the functor (\ref{top-component}).}
\label{zhu-left-adj-section}

Rosellen \cite{R} constructs the left adjoint
\begin{equation}
\label{zhu-left-adj} Zhu_V: Zhu(V)\rm{-}Mod \to V\rm{-}Mod
\end{equation}
to the  functor (\ref{top-component}) as follows.

% For details see [R], section 6.2.1.

\bigskip
To any vertex Lie algebra $R$ one can attach
  $$
  \fg(R)  = R[t,t^{-1}] / (\de a) (n) + n a (n-1),
  $$
a Lie algebra with bracket
  $$
  [a(n), b(m)] = \sum_{j\geq 0} {n \choose j} (a_{(j)} b )(n+m-j).
  $$
Here $a(n) = a t^n$, $a\in R$, $n\in \ZZ$.

\bigskip

If $V$ is a graded vertex algebra, $\fg(V)$ acquires a Lie algebra
grading with $a(n)$ sitting in component $n - \Delta_a+ 1$. Let us
concentrate on the subalgebra $\fg(V)_0$.

There is a surjective linear map $\fg(V)_0 \to Zhu(V)$ given by $
a_0 \mapsto [a] $.

\begin{lem}([R], Proposition 6.1.5)
  The map $a_0 \mapsto [a]$ from $\fg(V)_0 \to Zhu(V)$ is a Lie algebra map.
\end{lem}

%\bigskip

If $M$ is a $Zhu(V)$-module, then $M$ is a $\fg(V)_0$-module by
pullback (we have a natural Lie algebra map $\fg(V)_0 \to Zhu(V)$).

Moreover, we can extend this to a $\fg(V)_{\geq }$-module structure
by setting $\fg(V)_{>}M =\nolinebreak 0$. Then
\begin{equation}
  \tilde{M}  = U(\fg(V)) \otimes_{U(\fg(V)_{\geq})} M
\end{equation}
is a  $\ZZ_{\geq 0}$-graded   $\fg(V)$-module.

Let $Q(\tilde{M}) $  be the $\fg(V)$-submodule of $\tilde{M}$ generated by % the
 coefficients of
 \begin{center}
$ (a\opm  b)(z) m - :\nolinebreak a(z)b(z):\nolinebreak m,
  \ \  m\in \tilde{M} / Q(\tilde{M}). $
 \end{center}

Then $\tilde{M} / Q(\tilde{M})$ is a $\ZZ_{\geq 0}$-graded
$V$-module.

By definition,
\begin{equation}
  Zhu_V(M) = \tilde{M} / Q(\tilde{M})
\end{equation}
%%% In section 4 we will define sheaf versions of the functors above

\section{Universal  twisted cdo}
\subsection{Truncated de Rham complexes}\label{truncatedderham}\nopagebreak
Let $X$ be a smooth algebraic variety. For $0\leq p<q\leq
\text{dim}X$ introduce the complexes
$$
\Omega_X^{[p,q>}:\;
0\rightarrow\Omega^p_X\rightarrow\Omega^{p+1}_X\rightarrow\cdots\Omega^{q-1}_X
\rightarrow\Omega_X^{q,cl}\rightarrow 0,
$$
$$
\Omega_X^{[p}:\;
0\rightarrow\Omega^p_X\rightarrow\Omega^{p+1}_X\rightarrow\cdots,
$$
$$
\widetilde{\Omega}_X^{[q}:\;
0\rightarrow\frac{\Omega^q_X}{\Omega_X^{q,cl}}\rightarrow\Omega^{q+1}_X\rightarrow\cdots,
$$
where $\Omega_X^m$ stands for the sheaf of $m$-forms, the
differential is assumed to be the de Rham differential, and the
grading is shifted so that $\Omega_X^p$ is placed in degree 0.

For any complex of sheaves $\cA$ over $X$ consider the
hypercohomology groups $H^{i}(X,\cA)$, for any cover of $X$,
$\frak{U}=\{U_i\}$, the corresponding Cech hypercohomology,
$\check{H}^i(\frak{U},\cA)$, and finally the Cech hypercohomology
$\check{H}^i(X,\cA) =\lim_{\rightarrow}\check{H}^i(\frak{U},\cA)$.

The diagram
$$
\Omega^{[p,q>}_X\rightarrow\Omega^{[p}\rightarrow
\widetilde{\Omega}_X^{[q}
$$
is an exact triangle. The corresponding long cohomology sequence and
the fact that the de Rham cohomology
$H_{DR}^{m}(X,\widetilde{\Omega}_X^{[q})=0$ if $0\leq m\leq q-p$
implies, cf. \cite{GMSII}, sect.4.1.1,
\begin{lem}
\label{fact-on-coho}
 The canonical maps
$$
H^i(X,\Omega^{[p,q>})\rightarrow H^i(X,\Omega^{[p}),\;
\check{H}^i(X,\Omega^{[p,q>})\rightarrow \check{H}^i(X,\Omega^{[p})
$$
are isomorphisms if $0\leq i\leq q-p$ and  injections if $i=q-p+1$.
\end{lem}
\begin{cor}
\label{cech-derived} The canonical map
$$
\check{H}^i(X,\Omega^{[p,q>}_X)\rightarrow H^i(X,\Omega^{[p,q>}_X)
$$
is an isomorphism  if $i\leq q-p$.
\end{cor}
{\em Proof.} Since $\Omega_X^{[p}$ is an $\cO_X$-module,
$\check{H}^i(X,\Omega^{[p}_X)\rightarrow H^i(X,\Omega^{[p}_X)$ is an
isomorphism. The map indicated in the lemma factors as follows
$$
\check{H}^i(X,\Omega^{[p,q>}_X)\rightarrow
\check{H}^i(X,\Omega^{[p}_X)\rightarrow
H^i(X,\Omega^{[p}_X)\rightarrow H^i(X,\Omega^{[p,q>}_X),
$$
where thanks to Lemma \ref{fact-on-coho} all maps are isomorphisms
if $i\leq q-p$. $\qed$

\subsection{Twisted differential operators.}
\label{Twisted_differential_operators} A sheaf of twisted
differential operators (TDO) is a sheaf of filtered $\cO_X$-algebras
such that the corresponding graded sheaf is (the push-forward of)
$\cO_{T^*X}$, see [BB2]. The set of isomorphism classes of such
sheaves is in 1-1 correspondence with $H^1(X,\Omega^{[1,2>}_X)$.
Denote by $\cD^\la_X$ a TDO that corresponds to $\la\in
H^1(X,\Omega^{[1,2>}_X)$. If $\text{dim}
H^1(X,\Omega^{[1,2>}_X)<\infty$, then it is easy to construct a
universal TDO, that is to say, a family of sheaves with base
$H^1(X,\Omega^{[1,2>}_X)$ so that the sheaf that corresponds to a
point $\la\in H^1(X,\Omega^{[1,2>}_X)$ is isomorphic to $\cD^\la_X$.
The construction is as follows.

{\em Assume that $X$ is projective}. Then, as Lemma
\ref{fact-on-coho} implies,
$\text{dim}H^1(X,\Omega^{[1,2>})<\infty$. According to Corollary
\ref{cech-derived}, we can pick an affine cover $\frak{U}$ so that
$\check{H}^1(\frak{U},\Omega^{[1,2>})=H^1(X,\Omega^{[1,2>})$.

Let $\{ \la_k \}$  and $\{ \la^*_k \}$ be dual bases of
$\check{H}^1(\frak{U},\Omega^{[1,2>})$ and
$(\check{H}^1(\frak{U},\Omega^{[1,2>}))^*$ respectively. We fix a
lifting $\sigma: H^1(\frak{U}, \Omega^{[1,2>}) \to Z^1(\frak{U},
\Omega^{[1,2>})$ and identify the former with a subspace of the
latter using this lifting. Upon this identification, each $\la_k$
becomes a pair
\begin{equation}
 \label{definitionofacocycle}
\la_k=(\la_k^{(1)},\la_k^{(2)})\in(\prod_{i,j}\Gamma(U_i\cap
U_j,\Omega^1_X))\times (\prod_{i}\Gamma(U_i,\Omega^{2,cl}_X))
\end{equation}
so that the forms $\la_k^{(1)}(U_i\cap U_j)\in \Gamma(U_i\cap
U_j,\Omega^1_X)$ and $\la_k^{(2)}(U_{i})\in
\Gamma(U_i,\Omega^{2,cl}_X)$ are defined for each $k$, $i$, and $j$.
The cocycle condition reads
\begin{equation}
\label{cocycle_condition}
d_{DR}\la_k^{(1)}=d_{\check{C}}\la^{(2)}_k,\;
d_{\check{C}}\la_k^{(1)}=0.
\end{equation}

The space $\cO_{U_i}\otimes  H^1(X,\Omega^{[1,2>})^*$ carries two obvious
 actions, by $\cO_{U_i}$ and $\cT_{U_i}$, defined as follows
 \[\cO_{U_i}\otimes(\cO_{U_i}\otimes  H^1(X,\Omega^{[1,2>})^*)\rightarrow
 \cO_{U_i}\otimes  H^1(X,\Omega^{[1,2>})^*,\; f\otimes(g\otimes h)\mapsto
 fg\otimes h,
 \]
\[\cT_{U_i}\otimes(\cO_{U_i}\otimes  H^1(X,\Omega^{[1,2>})^*)\rightarrow
 \cO_{U_i}\otimes  H^1(X,\Omega^{[1,2>})^*,\; \tau\otimes(g\otimes h)\mapsto
 \tau(g)\otimes h,
 \]
and the actions are compatible in that
 $\tau(f\cdot p)=\tau(f)\cdot p+f\cdot\tau(p)$, $\tau\in\cT_{U_i}$, $f\in\cO_{U_i}$,
 $p\in\cO_{U_i}\otimes  H^1(X,\Omega^{[1,2>})^*$.

 Consider an abelian extension of
 $\cT_{U_i}$ by $\cO_{U_i}\otimes  H^1(X,\Omega^{[1,2>})^*$
 \begin{equation}
 \label{def-of-tau-tw-1}
 0\rightarrow \cO_{U_i}\otimes
 H^1(X,\Omega^{[1,2>})^*\rightarrow\cT^{tw}_{U_i}\rightarrow\cT_{U_i}\rightarrow 0
 \end{equation}
defined by the following cocycle
$$
\cT_{U_i}\ni\xi,\eta\mapsto \sum_k \iota_{\xi}\iota_{\eta}\la_k^{(2)}(U_i))\la^*_k.
$$
In other words, let us define the bracket $\widetilde{[.,.]}$ so
that
\begin{equation}
\label{deform_bracket} \widetilde{[\xi,\eta]}=
[\xi,\eta]+\sum_k(\iota_{\xi}\iota_{\eta}\la_k^{(2)}(U_i))\la^*_k,
\end{equation}
for all  $\xi,\eta\in \Gamma(U_i,\cT_{U_i})$; here $[\xi,\eta]$ is
the usual Lie bracket of vector fields.

The fact that each $\lambda^{(2)}_k(U_i)$ is a closed 2-form implies that $\cT^{tw}_{U_i}$ is
indeed a Lie algebra,  in fact a $\cO_{U_i}$-Lie algebroid,
$\cT^{tw}_{U_i}\rightarrow\cT_{U_i}$ being the anchor map. Let
$$
\cD^{tw}_{U_i}=U_{ {\cO}_{U_i} }    \cT^{tw}_{U_i},
$$
where $ U_{ {\cO}_{U_i} }$ is the enveloping algebra
functor, cf. the end of sect.\ \ref{Chiral_differential_operators}.

Define the transition maps $g_{ij}: \cD^{tw}_{U_i}|_{U_i\cap
U_j}\rightarrow \cD^{tw}_{U_j}|_{U_i\cap U_j}$ by requiring that
\begin{equation}
\label{deform_gluing} g_{ij}(\xi) = \xi -
\sum_k(\iota_{\xi}\la_k^{(1)}(U_i\cap U_j))\la^*_k, \ \ g_{ij}(f)=
f, \ \ f\in \cO_{U_i}\otimes\pole[H^1(X,\Omega^{[1,2>})].
\end{equation}
The  condition $d_{DR}\la_k^{(1)}=d_{\check{C}}\la^{(2)}_k$ implies
that each $g_{ij}$ is an associative algebra homomorphism, and the
condition $d_{\check{C}}\la_k^{(1)}=0$ implies that
$g_{ij}=g_{ik}\circ g_{kj}$.  Denote by $\cD^{tw}_X$ the sheaf
obtained by gluing the sheaves $\cD_{U_i}^{tw}$ over two-fold
intersections via the  maps $g_{ij}$.

By construction,  $\pole[H^1(X,\Omega^{[1,2>})]$ lies in the center
of $\Gamma(X,\cD^{tw}_X)$, and if we let ${\frak m}_{\la}\in
\pole[H^1(X,\Omega^{[1,2>})]$ be the maximal ideal defined by
$\lambda\in H^1(X,\Omega^{[1,2>})$, then by definition,
\begin{equation}
\label{specializ-us-tw} \cD^{tw}_X/{\frak
m}_{\lambda}\cD^{tw}_X\text{ is isomorphic to }\cD^{\lambda}_X.
\end{equation}

It is clear that $\cD^{tw}_X$ is independent of the cover ${\frak
U}$ and lifting $\sigma$, and we call this sheaf  \textbf{the
universal sheaf of twisted differential operators}.

\subsection{Chiral analogue}

\subsubsection{A universal twisted CDO}
\label{A_universal_twisted_CDO}

Let $ch_2(X)=0$ and fix a CDO $\cD^{ch}_X$. To each such sheaf we
will attach a {\em universal twisted CDO}, $\cD^{ch, tw}_X$, in a
manner analogous to that in which we constructed a universal TDO
$\cD^{tw}_X$ in the previous section. Let us then place ourselves in
the situation of the previous section, where we had a fixed affine
cover ${\frak U}=\{U_i\}$ of a projective algebraic manifold $X$,
 dual bases $\{\lambda_i\}\in
H^1(X,\Omega^{[1,2>}_X)$, $\{\lambda_i^*\}\in
H^1(X,\Omega^{[1,2>}_X)^*$, and a lifting
$H^1(X,\Omega^{[1,2>}_X)\rightarrow Z^1({\frak
U},\Omega^{[1,2>}_X)$.

Assuming, as we may, that each $U_i$ is suitable for chiralization
we fix, for each $i$, an abelian basis
$\tau^{(i)}_1,\tau^{(i)}_2,...$ of $\Gamma(U_i,\cT_X)$, and a
collection of 3-forms $\alpha^{(i)}\in\Gamma(U_i,\Omega^{3,cl}_X)$,
cf. sect.\ \ref{Chiral_differential_operators},
Theorem~\ref{class-cdo-local}.
\begin{lem}
\label{local_model_chir}

(a) There is a unique vertex algebroid structure on the sheaf
\[
\cA_{U_i}\stackrel{\text{def}}{=}
\cO_{U_i}\oplus\cT_{U_i}\oplus\Omega_{U_i}\oplus\left(\oplus_{j}\cO_{U_i}\otimes\pole\lambda^*_j\right)
\]
so that
\[
\text{(1) } V_0=\cO_{U_i},\;
V_1=\cT_{U_i}\oplus\Omega_{U_i}\oplus\left(\oplus_{j}\cO_{U_i}\otimes\pole\lambda^*_j\right);
\]
\[
\text{(2) }\partial: \cO_{U_i}\rightarrow \Omega_{U_i}\text{ is the
de Rham differential};
\]
(3) the pair $(V_0,_{(-1)})$ is  $\cO_{U_i}$ as a commutative
associative algebra;
\[
\text{(4) } f_{(-1)}\omega=f\omega,\;
f\in\cO_{U_i},\omega\in\Omega_{U_i};
\]
\[
\text{(5) } f_{(-1)}\xi = f\xi\text{ mod
}\Omega_{U_i},\;f\in\cO_{U_i},\xi\in\cT_{U_i};
\]
\[
\text{(6)
}\tau^{(i)}_{l(0)}\tau^{(i)}_{m}=\iota_{\tau^{(i)}_{l}}\iota_{\tau^{(i)}_{m}}\alpha^{(i)}+
\sum_k(\iota_{\tau^{(i)}_{l}}\iota_{\tau^{(i)}_{m}}\lambda^{(2)}_k(U_i))
\lambda^*_k,\; \tau^{(i)}_{l(0)}f=\tau^{(i)}_l(f),f\in\cO_{U_i};
\]
\[
\text{(7) }\tau^{(i)}_{l(1)}\tau^{(i)}_{m}=0;
\]
\[
\text{(8) }\lambda^*_{k(0)}a=\lambda^*_{k  (1)}a = 0\text{ for any
}k, a.
\]

(b) The corresponding Lie algebroid $T=T(\cA_{U_i})$ satisfies,
\[
T=\cT_{U_i}^{tw},
\]
where $\cT_{U_i}^{tw}$ is the Lie algebroid that was defined in
sect.~\ref{Twisted_differential_operators}.
\end{lem}
{\em Proof.}

(a) It is clear, cf. sect.~\ref{Chiral_differential_operators}, that
there is only one way to extend the indicated operations
to the entire $\cA_{U_i}$  using the Borcherds identity (\ref{Borcherds-identity}).
Furthermore, thus obtained operations are all represented by
differential operators. In order to verify that these operations
satisfy the identities imposed by the definition of a vertex
algebroid, let us embed the sheaf in question, $\cA_{U_i}$, into its
formal completion, $\hat{\cA}_{U_i,x}$, at an arbitrary point $x\in
U_i$.  All operations on $\cA_{U_i}$ extend to those on $\hat{\cA}_{U_i,x}$.
 We will, first, prove that $\hat{\cA}_{U_i,x}$ with these operations is a
vertex  algebroid.

Upon passing to this completion each 2-form $\lambda^{(2)}_k(U_i)$ becomes
exact, and for each $k$ we obtain $\mu_k$ such that
$d_{DR}\mu_k=\lambda^{(2)}_k(U_i)$.  Now replace each $\tau^{(i)}_l$
with $\tilde{\tau}^{(i)}_l=
\tau^{(i)}_l+\sum_k\iota_{\tau^{(i)}_l}\mu_k\lambda^*_k$. It is
clear that in terms of this new basis condition (6) of our lemma
becomes
$$
\text{(6')
}\tilde{\tau}^{(i)}_{l(0)}\tilde{\tau}^{(i)}_{m}=\iota_{\tau^{(i)}_{l}}\iota_{\tau^{(i)}_{m}}\alpha^{(i)}.
$$
This means that the subspace spanned over $\Cplx$ by
$\lambda^{*}_j$, $1\leq j\leq n$, decouples. More precisely, if we
let
\[
\widetilde{\cA}_{U_i,x}=\hat{\cO}_{U_i,x}(\oplus_l
\pole \tilde{\tau}^{(i)}_l)\oplus\hat{\Omega}_{U_i,x},
\]
then the fact that $\tilde{\cA}_{U_i,x}$ with operations
(1--5,6',7,8) is a vertex algebroid becomes one of the main
observations of \cite{GMS}, recorded above as
Theorem~\ref{class-cdo-local}.

Adjoining the commutative variables $\lambda^*_j$ is easy. Condition
(8) above simply means that, as a space with operations $\partial$,
$_{(n)}$, $n=-1,0,1$,
\[
\hat{\cA}_{U_i,x}=\tilde{\cA}_{U_i,x}\stackrel{\bullet}{\otimes}(\pole\oplus
(\sum_j\pole\lambda^*_j)),
\]
where the tensor product functor is as in
(\ref{def-of-tens-prod-algebroid-vert}) and $\pole\oplus
(\sum_j\pole\lambda^*_j)$ is a commutative vertex algebroid from
sect. \ref{examples-of-vert-alg}. Since the R.H.S. is a vertex
algebroid, so is the L.H.S., $\hat{\cA}_{U_i,x}$.

The map $\cA_{U_i}\rightarrow\hat{\cA}_{U_i,x}$ being an injection,
the passage to the completion cannot create any new identities;
hence the operations initially defined on $\cA_{U_i}$ also satisfy
the definition of a vertex algebroid.

(b) It  was explained in sect.~\ref{def-of-vert-alg} that, as an
$\cO_{U_i}$-module,
\[
T=\left(\cT_{U_i}\oplus\Omega_{U_i}\oplus\left(\oplus_{j}\cO_{U_i}\otimes\pole\lambda^*_j\right)\right)/\Omega_{U_i},
\]
hence
\[
T=\cT_{U_i}\oplus\left(\oplus_{j}\cO_{U_i}\otimes\pole\lambda^*_j\right),
\]
which is precisely $\cT_{U_i}^{tw}$. The Lie bracket is defined by
$_{(0)}$. It remains to notice that upon  passing over to the
quotient modulo $\Omega_{U_i}$, formula (6) of Lemma
\ref{local_model_chir} becomes exactly formula
(\ref{deform_bracket}). $\qed$

\bigskip

Define a sheaf of vertex algebras over each $U_i$ by applying the
vertex enveloping algebra functor as follows
\begin{equation}
\label{finally-a-def-of-ctdo}
\cD^{ch,tw}_{U_i}\stackrel{\text{def}}{=}U\cA_{U_i}.
\end{equation}
These sheaves will serve as local models for the universal twisted
CDO we are after.

By construction we have sheaf embeddings
$$
\cO_{U_i}\oplus\cT_{U_i}\oplus(\oplus_j\pole\lambda^*_j)\hookrightarrow\cD^{ch,tw}_{U_i}.
$$
Recall now that we have assumed given a CDO $\cD^{ch}_X$. One way to
define this sheaf is to introduce the restrictions
$\cD^{ch}_{U_i}=\cD^{ch}_{X}|_{U_i}$, fix splittings
$$
\cO_{U_i}\oplus\cT_{U_i}\hookrightarrow\cD^{ch,tw}_{U_i},
$$
and the corresponding transition functions
$$
g_{ij}: \cD^{ch}_{U_i}|_{U_i\cap U_j}\rightarrow
\cD^{ch}_{U_j}|_{U_i\cap U_j}.
$$
\begin{lem}
 \label{gluing_chir}

(1) There is a unique vertex algebra isomorphism
$$
g_{ij}^{tw}: \cD^{ch,tw}_{U_i}|_{U_i\cap U_j}\rightarrow
\cD^{ch,tw}_{U_j}|_{U_i\cap U_j}
$$
such that
$$
g_{ij}^{tw}|_{\cO_{U_i\cap U_j}}=g_{ij}|_{\cO_{U_i\cap U_j}}
$$
$$
g_{ij}^{tw}(\lambda_k^*)=\lambda_k^*,
$$
$$
g_{ij}^{tw}(\xi)=g_{ij}(\xi)-\sum_k(\iota_{\xi}\lambda^{(1)}_{k}(U_{ij}))\lambda_k^*,
\xi\in\cT_{U_i\cap U_j}.
$$
(2) On triple intersections $U_i\cap U_j\cap U_k$
$$
g^{tw}_{ij}=g^{tw}_{kj}\circ g^{tw}_{ik}.
$$
\end{lem}

\bigskip

{\em Proof.}

(1)As it follows from the {\em Reconstruction theorem}, [K], the
fact on which an analogous discussion in \cite{MSV} heavily relies,
it is enough to verify the equalities
$$
g_{ij}^{tw}(\xi)_{(1)}g_{ij}^{tw}(\eta)=g_{ij}^{tw}(\xi_{(1)}\eta),\;
g_{ij}^{tw}(\xi)_{(0)}g_{ij}^{tw}(\eta)=g_{ij}^{tw}(\xi_{(0)}\eta),\;\xi,\eta\in\cT_{U_i\cap
U_j}.
$$
The former is part of the definition of $\cD^{ch}_X$ for we have, by
definition,
$g_{ij}^{tw}(\xi)_{(1)}g_{ij}^{tw}(\eta)=g_{ij}(\xi)_{(1)}g_{ij}(\eta)$
and $g_{ij}^{tw}(\xi_{(1)}\eta)=g_{ij}(\xi_{(1)}\eta)$.

The latter  boils down to the purely classical statement that
underlies the construction of the twisted differential operators,
see sect.\ \ref{Twisted_differential_operators}.  Note that the
deformation of $_{(0)}$ by a function introduced in Lemma
\ref{local_model_chir}(6) has as a consequence the fact that the
``old'' transition functions, $g_{ij}$, are no longer vertex algebra
morphisms, the discrepancy being
$$
g_{ij}(\xi_{(0)}\eta)-g_{ij}(\xi)_{(0)}g_{ij}(\eta)=\sum_k\iota_{\xi}\iota_{\eta}(\lambda^{(2)}_k(U_i)-
\lambda^{(2)}_k(U_j))\lambda^*_k.
$$
This discrepancy is taken care of by the passage from $g_{ij}$ to
$g_{ij}^{tw}$. Indeed, since by definition
$$
g_{ij}(\xi)_{(0)}\sum_k(\iota_{\eta}\lambda^{(1)}_{k}(U_{ij}))\lambda_k^*=
-(\sum_k(\iota_{\eta}\lambda^{(1)}_{k}(U_{ij}))\lambda_k^*)_{(0)}g_{ij}(\xi)=
\sum_k\xi(\iota_{\eta}\lambda^{(1)}_{k}(U_{ij}))\lambda_k^*,
$$
we have
$$
g_{ij}^{tw}(\xi_{(0)}\eta)=g_{ij}(\xi_{(0)}\eta)-\sum_k\iota_{[\xi,\eta]}\lambda^{(1)}_k(U_{ij})\lambda^*_k;
$$
$$
g_{ij}^{tw}(\xi)_{(0)}g_{ij}^{tw}(\eta)=
g_{ij}(\xi)_{(0)}g_{ij}(\eta)-\sum_k\xi((\iota_{\eta}\lambda^{(1)}_{k}(U_{ij}))\lambda_k^*+
\sum_k\eta((\iota_{\xi}\lambda^{(1)}_{k}(U_{ij}))\lambda_k^*.
$$
Subtracting the latter from the former we obtain
$$
g_{ij}^{tw}(\xi_{(0)}\eta)-g_{ij}^{tw}(\xi)_{(0)}g_{ij}^{tw}(\eta)=
$$
$$
\sum_k ( \iota_{\xi}\iota_{\eta}(\lambda^{(2)}_k(U_i)-\lambda^{(2)}_k(U_j)) +
\sum_k  (\iota_{\eta}\iota_{\xi}  d_{DR}\lambda^{(1)}_k(U_{ij}))\lambda^*_k=
$$
$$
-\sum_k(d_{\check{C}}(\lambda^{(2)}_k(U_{ij}))|_{\xi,\eta}-d_{DR}\lambda^{(1)}_k(U_{ij})|_{\xi,\eta})\lambda^*_k
$$
which vanishes by virtue of the first part of the cocycle condition
(\ref{cocycle_condition}).

(2) This is also a statement about twisted differential operators:
we have over $U_i\cap U_j\cap U_k$
$$
g^{tw}_{ij}(\xi)-g^{tw}_{kj}\circ
g^{tw}_{ik}(\xi)=g_{ij}(\xi)-g_{kj}\circ g_{ik}(\xi)-
\sum_k(d_{\check{C}}\lambda^{(1)}_k|_{\xi})\lambda^*_k,
$$
which vanishes by virtue of the second part of the cocycle condition
(\ref{cocycle_condition}). $\qed$

\bigskip

Lemma \ref{gluing_chir} implies the following
\begin{thm-def}
\label{thm-def-twist-chir} Given a projective algebraic manifold $X$
and a CDO $\cD^{ch}_X$,
 there is a unique sheaf of vertex algebras, to be denoted $\cD^{ch,tw}_X$ and called
{\bf a universal sheaf of twisted chiral differential operators}
(TCDO), such that
$$
\cD^{ch,tw}_X|_{U_i}=\cD^{ch,tw}_{U_i},
$$
and the canonical isomorphisms
$$
\cD^{ch,tw}_{U_i}|_{U_i\cap U_j}\rightarrow
\cD^{ch,tw}_{U_j}|_{U_i\cap U_j}
$$
are $g_{ij}^{tw}$.
\end{thm-def}

\bigskip
Indeed, the assumptions of the theorem require that $\cD^{ch,tw}_X$
be obtained by gluing the pieces $\cD^{ch,tw}_{U_i}$ via
$g_{ij}^{tw}$, and the gluing is made sense of by Lemma
\ref{gluing_chir}. $\qed$

\bigskip

Let $H_X$ be the commutative vertex algebra of differential
polynomials on $H^1(X,\Omega^{[1,2>}_X)$, cf. sect.\
\ref{Commutative_vertex_algebras}. As a commutative algebra
$$
H_X=\pole[\partial^j\lambda^*_k; j\geq 0,1\leq k\leq
\text{dim}H^1(X,\Omega^{[1,2>}_X)],
$$
and the canonical derivation is $\partial$.

Denote by $\uH_X$  the constant sheaf over $X$ with $\uH_X(U) = H_X$
for nonempty $U$.

It is clear that  if we let operations $_{(0)}, _{(1)}=0$, then
$\pole\oplus(\oplus_j\pole\lambda^*_j)$ is a vertex algebroid  and
that $U(\pole\oplus(\oplus_j\pole\lambda^*_j)=H_X$. Now the
embeddings
\begin{equation}
 \label{center}
\uH_X\hookrightarrow Z(\cD^{ch,tw}_X),\; H_X\hookrightarrow
Z(\Gamma(X,\cD^{ch,tw}_X))
\end{equation}
follow from the constructions at once; here for any vertex algebra
$V$, $Z(V)$ stands for its center, that is to say, $Z(V)=\{v\in
V\text{ s.t. }v_{(n)}V=0\text{ for all }n\geq 0\}$.

\subsubsection{Locally trivial and other versions of twisted CDOs}
\label{Locally_trivial_and_other_versions_of_twisted_CDOs}

To begin with,  note that the requirement that $X$ be projective was
needed above only to ensure that $H^1(X,\Omega^{[1,2>}_X)$ is
finite-dimensional. In the infinite-dimensional situation one has to
work with completions, which may well be possible but not
attractive.

On the other hand, for any $X$ and a fixed cover ${\frak U}$ one can
repeat all of the above constructions and obtain sheaves
$\cD^{tw}_{X,{\frak U}}$ and $\cD^{ch,tw}_{X,{\frak U}}$. Such
sheaves will not be universal in general and will explicitly depend
on the choice of ${\frak U}$.

Yet another version of our construction will give us locally trivial
twisted sheaves of chiral differential operators.

There is an embedding
\begin{equation}
\label{emb_coho}
 H^{1}(X,\Omega^{1,cl}_X)\hookrightarrow H^{1}(X,\Omega^{[1,2>}_X)
\end{equation}
The space $H^{1}(X,\Omega^{1,cl}_X)$ classifies {\em locally
trivial} twisted differential operators, those that are locally
isomorphic to $\cD_X$. Thus for each $\lambda\in
H^{1}(X,\Omega^{1,cl}_X)$, there is a unique up to isomorphism TDO
$\stackrel{\circ}{\cD}^{\lambda}_X$ such that for each sufficiently
small open $U\subset X$, $\stackrel{\circ}{\cD}^{\lambda}_X|_U$ is
isomorphic to $\cD_U$. Let us  see what this means at the level of
the universal TDO.

In terms of Cech cocycles the image of embedding (\ref{emb_coho}) is
described by those $(\lambda^{(1)},\lambda^{(2)})$, see
(\ref{definitionofacocycle}), where $\lambda^{(2)}=0$, and this
forces $\lambda^{(1)}$ to be closed.  Picking a collection of such
cocycles that represent a basis of $H^{1}(X,\Omega^{1,cl}_X)$ we can
repeat the constructions of sections
\ref{Twisted_differential_operators} and
\ref{A_universal_twisted_CDO} to obtain sheaves
$\stackrel{\circ}{\cD}^{tw}_X$ and
$\stackrel{\circ}{\cD}^{ch,tw}_X$. The former is glued of the pieces
isomorphic to $\cD_{U_i}\otimes\pole[H^{1}(X,\Omega^{1,cl}_X)]$ as
associative algebras (this is a locally trivial property, it is due
to the vanishing of $\lambda^{(2)}$), the transition functions being
defined as in (\ref{deform_gluing}). The latter is defined likewise
by gluing pieces isomorphic (as vertex algebras) to
$\cD^{ch}_{U_i}\otimes H_X$ with transition functions as in Lemma
\ref{gluing_chir}; here $H_X$ is the vertex algebra of differential
polynomials on
 $H^{1}(X,\Omega^{1,cl}_X)$. We will call the sheaves $\stackrel{\circ}{\cD}^{tw}_X$
  and $\stackrel{\circ}{\cD}^{ch,tw}_X$ the {\em universal locally trivial sheaves of twisted (chiral resp.) differential
operators.}
%
%If ${\frak m}_{\lambda}$ is the maximal ideal of $\pole[ H^{1}(X,\Omega^{1,cl}_X)]$ corresponding to the point $\lambda\in H^{1}(X,\Omega^{1,cl}_X)$, then one obtains the following analogue of
%(\ref{specializ-us-tw})
%\begin{equation}
% \label{specializ-us-tw-loctr}
%\stackrel{\circ}{\cD}^{tw}_X/{\frak m}_{\lambda}\stackrel{\circ}{\cD}^{tw}_X.
%\end{equation}

\subsubsection{Example: flag manifolds.}
\label{Example:_flag_manifolds.} Let us see what our constructions
give us if $X=\prline$. We have
$\prline=\pole_{0}\cup\pole_{\infty}$, a cover ${\frak
U}=\{\pole_{0},\pole_{\infty}\}$, where $\pole_0$ is $\pole$ with
coordinate $x$, $\pole_{\infty}$ is $\pole$ with coordinate $y$,
with the transition function $x\mapsto 1/y$ over
$\pole^*=\pole_{0}\cap\pole_{\infty}$.

Defined over $\pole_{0}$ and $\pole_{\infty}$ are the standard CDOs,
$\cD^{ch}_{\pole_0}$ and $\cD^{ch}_{\pole_{\infty}}$. The spaces of
global sections of these sheaves are  polynomials in
$\partial^{n}(x)$, $\partial^{n}(\partial_x)$ (or $\partial^{n}(y)$,
$\partial^{n}(\partial_y)$ in the latter case), where $\partial$ is
the translation operator, so that, cf. sect.\
\ref{Chiral_differential_operators},
$$
(\partial_x)_{(0)}x=(\partial_y)_{(0)}y=1.
$$
There is a unique up to isomorphism CDO on $\prline$,
$\cD^{ch}_{\prline}$; it is defined by gluing $\cD^{ch}_{\pole_0}$
and $\cD^{ch}_{\pole_{\infty}}$ over $\pole^*$ as follows
\cite{MSV}:
\begin{equation}
 \label{gluing_untwisted_p1}
x\mapsto 1/y,\; \partial_x\mapsto
(-\partial_{y})_{(-1)}(y^2)-2\partial(x).
\end{equation}
The canonical Lie algebra morphism
\begin{equation}
\label{sl2-p1} sl_2\rightarrow \Gamma(\prline,\cT_{\prline}),
\end{equation}
where
\begin{equation}
 \label{formulas-sl2-p1}
 e\mapsto\partial_x,  \quad   h\mapsto -2x\partial_x,  \quad  f\mapsto -x^2\partial_x,
\end{equation}
$e,h,f$ being the standard generators of $sl_2$, can be lifted to a
vertex algebra morphism
\begin{equation}
\label{verte-sl2-p1} V_{-2}(sl_2)\rightarrow
\Gamma(\prline,\cD^{ch}_{\prline}),
\end{equation}
where
\begin{equation}
 \label{verte-formulas-sl2-p1}
e_{(-1)}\vac\mapsto\partial_x,
 h{(-1)}\vac\mapsto -2(\partial_x)_{(-1)}x,
 f_{(-1)}\vac\mapsto -(\partial_x)_{(-1)}x^2-2\partial(x).
\end{equation}
The  twisted version of all of this is as follows.

Since $\text{dim}\prline=1$,
\[
H^1(\prline,\Omega^{1}_{\prline}\rightarrow\Omega^{2,cl}_{\prline})=
\Omega^{1}_{\prline}=\Omega^{1,cl}_{\prline},
\]
so all twisted CDO on $\prline$ are locally trivial. Furthermore,
$H^1(\prline,\Omega^{1,cl}_{\prline})=\pole$ and is spanned by the cocycle
$\pole_0\cap\pole_{\infty} \mapsto dx/x$, the Chern class of Serre's sheaf $\cO(1)$. We have
$H_{\prline}=\pole[\lambda^*,\partial(\lambda^*),....]$. Let
$\cD^{ch,tw}_{\pole_0}=\cD^{ch}_{\pole_0}\otimes H_{\prline}$,
$\cD^{ch,tw}_{\pole_{\infty}}=\cD^{ch}_{\pole_{\infty}}\otimes H_{\prline}$ and define $
\cD^{ch,tw}_{\prline}$ by gluing $\cD^{ch,tw}_{\pole_0}$ onto $ \cD^{ch,tw}_{\pole_{\infty}}$
via
\begin{equation}
 \label{gluing_twisted_p1}
\lambda^*\mapsto\lambda^*,\;x\mapsto 1/y,\;\partial_x\mapsto
-(\partial_y)_{(-1)}y^2-2\partial(y) +y_{(-1)}\lambda^*.
\end{equation}
Morphism (\ref{verte-sl2-p1}) ``deforms'' to
\begin{equation}
\label{tw-verte-sl2-p1} V_{-2}(sl_2)\rightarrow \Gamma(\prline, \cD^{ch,tw}_{\prline}),
\end{equation}
\begin{equation}
 \label{tw-verte-formulas-sl2-p1}
e_{(-1)}\vac\mapsto\partial_x, h{(-1)}\vac\mapsto
-2(\partial_x)_{(-1)}x+\lambda^*, f_{(-1)}\vac\mapsto
-(\partial_x)_{(-1)}x^2-2\partial(x)+x_{(-1)}\lambda^*.
\end{equation}
Furthermore, consider
$T=e_{(-1)}f_{(-1)}+f_{(-1)}e_{(-1)}+1/2h_{(-1)}h\in V_{-2}(sl_2)$.
It is known that $T\in\fz(V_{-2}(sl_2))$, the center of
$V_{-2}(sl_2)$, and in fact, the center $\fz(V_{-2}(sl_2))$ equals
the commutative vertex algebra of differential polynomials in $T$.
The formulas above show
\begin{equation}
\label{image-sugawara}
T\mapsto\frac{1}{2}\lambda^*_{(-1)}\lambda^*-\partial(\lambda^*)\in
H_{\prline}.
\end{equation}

All of the above is easily verified by direct computations, cf. \cite{MSV}. The higher rank
analogue is less explicit but valid nevertheless.

\bigskip

Let $G$ be a simple complex Lie group, $B\subset G$ a Borel subgroup, $X=G/B$, the flag
manifold,  $\fg=\text{Lie\;}G$ the corresponding Lie algebra, $\fh$ a Cartan subalgebra.  One
has a sequence of maps
\begin{equation}
 \fh^* \rightarrow H^1(X,\Omega^{1,cl}_X)\rightarrow
H^1(X,\Omega^{1}_X\rightarrow\Omega^{2,cl}_X).
\end{equation}
The leftmost map    attaches to an integral weight $\lambda\in P\subset\fh^*$ the Chern class
of the $G$-equivariant line bundle $\cL_\lambda=G\times_{B}\pole_{\lambda}$, and then extends
thus defined map $P\rightarrow H^1(X,\Omega^{1,cl}_X)$ to $\fh^*$ by linearity. The rightmost
one is engendered by the standard spectral sequence converging to  hypercohomology. It is easy
to verify that both these maps are isomorphisms. Therefore,
\begin{equation}
\label{on-flag-trivial} \fh^* \stackrel{\sim}{\rightarrow}H^1(X,\Omega^{1,cl}_X)
\stackrel{\sim}{\rightarrow}H^1(X,\Omega^{1}_X\rightarrow\Omega^{2,cl}_X),
\end{equation}
and each twisted CDO on $X$ is locally trivial.

Note that $\cL_\lambda$ being $G$-equivariant, there arises a map from $U\fg$ to the algebra of
differential operators acting on $\cL_\lambda$ or, equivalently, \cite{BB2},
\[
U\fg\rightarrow \cD_X^\lambda.
\]
A moment's thought shows that this map is a polynomial in $\lambda$;
hence it defines a universal map
\begin{equation}
\label{morph-bb-twisted} U\fg\rightarrow \cD_X^{tw}.
\end{equation}
 Constructed in \cite{MSV} is a (unique up to isomorphism \cite{GMSII}) CDO
$\cD^{ch}_X$. We arrive at the universal  twisted CDO $\cD^{ch,tw}_X$ locally isomorphic to
$\cD^{ch}_U\otimes H_X$, where $H_X$ is the vertex algebra of differential polynomials on
$\fh^*$.

Constructed in \cite{MSV} -- or rather in \cite{FF1}, see also \cite{F1} and \cite{GMSII} for
an alternative approach -- is a vertex algebra morphism
\begin{equation}
 \label{verte-g-f}
V_{-h^{\vee}}(\fg)\rightarrow\Gamma(X,\cD^{ch}_X).
\end{equation}
Furthermore, it is an important result of Feigin and Frenkel \cite{FF2}, see also an excellent
presentation in \cite{F1}, that $V_{-h^{\vee}}(\fg)$ possesses a non-trivial center,
$\fz(V_{-h^{\vee}}(\fg))$, which, as a vertex algebra, isomorphic to the algebra of
differential polynomials in $\text{rk}\fg$ variables.
\begin{lem}
 \label{tw-verte-g-f}
Morphism (\ref{verte-g-f}) ``deforms'' to
$$
\rho: \;V_{-h\check{ }}(\fg)\rightarrow\Gamma(X,\cD^{ch,tw}_X)
$$
and
$$
\rho(\fz(V_{-h^{\vee}}(\fg)))\subset H_{X}.
$$
\end{lem}

\bigskip
{\it Sketch of Proof.} We will be brief, because this is one of
those proofs that the reader may find easier to find on his own than
to read somebody else's explanations. For each $x\in\fg$, $\rho(x)$
can be written, schematically, as follows
$$
\rho(x)=(classical)+(chiral)+(classical)_{\lambda},
$$
where $(classical)$ are those terms that appear in the image of the
canonical map $U\fg\rightarrow \cD_{X}$,
$(classical)+(classical)_{\lambda}$ are those  that appear in the
image of the Beilinson-Bernstein map (\ref{morph-bb-twisted}), and
$(chiral)$ is the rest; note that equivalently
$(classical)+(chiral)$ is the image of map (\ref{verte-g-f}).

We have to verify that $\rho(x)_{(1)}\rho(y)=-h^{\vee}<x,y>$ and
$\rho(x)_{(0)}\rho(y)=\rho([x,y])$. Only terms
$(classical)+(chiral)$ contribute to the former; that their
contribution is as needed is the content of assertion
(\ref{verte-g-f}). Given the former, the latter becomes precisely
the classical construction of the morphism $U\fg\rightarrow
\cD_{X}^{tw}$.

\begin{sloppypar}
The assertion on the image of the center was actually verified in \cite{FF2,F1}. Indeed, since
$H_{X}$ is the space of global sections of the constant sheaf $\uH_X$, it is enough to verify
the assertion for the composition of $\rho$ with the embedding of $\Gamma(X,\cD^{ch,tw}_X)$ in
$\Gamma(X_e,\cD^{ch,tw}_X)$, where $X_e\subset X$ is the big cell. The space\newline
$\Gamma(X_e,\cD^{ch,tw}_X)$ is a Wakimoto module, and it is the properties of thus defined
morphism from $V_{-h\check{ }}(\fg)$ to the Wakimoto module that were studied in \cite{FF2,F1}.
$\qed$
\end{sloppypar}

\bigskip

\subsection{The Zhu algebra of $\twcdo$.}

Now we compute the Zhu algebra for the sheaf $\twcdo$. We show that
the obtained sheaf is the universal sheaf $\gtdo$ of twisted
differential operators on $X$.

\begin{thm}
\label{zhutcdo} Suppose $\cdo$  is a  CDO on $X$.
 Let $\twcdo$ be the corresponding twisted sheaf.
Then
 \begin{equation}
 Zhu(\twcdo)= \cD^{tw}_X.
 \end{equation}
 Likewise
 \begin{equation}
 Zhu(\twcdoloc)= \stackrel{\circ}{\cD}^{tw}_X.
 \end{equation}
\end{thm}

\textbf{Proof}. Let us compute $Zhu(\twcdo(U))$ for  $U_i\in
{\frak{U}}$.

 By definition, see (\ref{finally-a-def-of-ctdo}), we have $\cD^{ch,tw}_{U_i} = U\cA_{U_i}$.

Lemma \ref{local_model_chir} (b) says that the corresponding Lie
algebroid is $\cT^{tw}_{U_i}$.  Now Theorem \ref{zhu-envelope}
implies that $Zhu(\cD^{ch,tw}_{U_i}) = U_{\cO_{U_i}}\cT^{tw}_{U_i}$.
The latter by definition equals $\cD^{tw}_{U_i}$.

It remains to show that the transition functions are as claimed, and
this is obvious.

Literally the same proof applies to the locally trivial TCDO
$\twcdoloc$. $\qed$

\section{Modules over a universal twisted CDO}
\label{Modules_over_universal_twisted_CDO}

\subsection{The main result}
\label{The_main_result} We will call a sheaf of vector spaces $\cM$
a $\cD^{ch,tw}_X$-{\em module} if

(1) for each open $U\subset X$ is a
$\Gamma(U,\cD^{ch,tw}_X)$-module;

(2) the restriction morphisms
$\Gamma(U,\cM)\rightarrow\Gamma(V,\cM)$, $V\subset U$, are
 are  $\Gamma(U,\cD^{ch,tw}_X)$-module morphisms, where the
 $\Gamma(U,\cD^{ch,tw}_X)$-module structure on $\Gamma(V,\cM)$ is that of the pull-back w.r.t. to
the restriction map
$\Gamma(U,\cD^{tw}_X)\rightarrow\Gamma(V,\cD^{tw}_X)$;

 (3) $\cM$ is generated by a subsheaf
 $\cM_0$
 such that for each open $U\subset X$
 \begin{equation}
 \label{generated_by_Mzero}
 v_n\Gamma(U,\cM_0) = 0 \ \ \ \textrm{for } v\in \Gamma(U,\cD^{ch,tw}_X), \
 n>0,
 \end{equation}
 \begin{equation}
 \label{preserved_by_Mzero}
 v_0\Gamma(U,\cM_0) \subset \Gamma(U,\cM_0)\ \ \ \textrm{for } v\in
 \Gamma(U,\cD^{ch,tw}_X).
 \end{equation}

\begin{rem}
  Note that condition (3) implies  a $\twcdo$-module $\cM$ is  filtered, i.e.
%\end{rem}
\begin{equation}
\label{filtr-module-1}
\{0\}\subset\cM_0\subset\cM_1\subset\cdots,\;\cM=\cup_{n=0}^{\infty}\cM_n,\text{
with } \cM_j\stackrel{\text{def}}{=}\sum_{i=0}^{j}(\twcdo)_i\cM_0,
\end{equation}
and this filtration is compatible with the conformal weight grading
of $\cD^{ch,tw}_X$ in that
\begin{equation}
 \label{filtr-module-2}
((\cD^{ch,tw}_X)_{j})_{(l)}\cM_n\subset\cM_{j+n-l-1}.
\end{equation}
\end{rem}
\bigskip
Denote by $Mod-\cD^{ch,tw}_X$ the category of
$\cD^{ch,tw}_X$-modules.

Precisely the same definition can be made in the case of a locally
trivial TCDO, $\twcdoloc$, see sect.\
\ref{Locally_trivial_and_other_versions_of_twisted_CDOs} and we
obtain the category $Mod-\twcdoloc$.

Recall that $\cD^{ch,tw}_X$ contains a huge center $H_X\subset
Z(\Gamma(X,\cD^{ch,tw}_X))$, see (\ref{center}). Since the vertex
algebra $H_X$ is commutative, its irreducibles are all
1-dimensional, characters in other  words, and are in 1-1
correspondence with the algebra of Laurent series with values in
$H^1(X,\Omega^{[1,2>}_X)$. Specifically, if  $\chi(z)\in
H^1(X,\Omega^{[1,2>}_X)((z))$, then the character $\pole_{\chi}$ is
a 1-dimensional $H_X$-module defined by, cf (\ref{def-vert-mod-1},
\ref{def-vert-mod-2})
\begin{equation}
\label{def-of-char} \chi: H_X\rightarrow
Fields(\Cplx_\chi),\;\chi(\lambda)(z)=\lambda(\chi(z)).
\end{equation}
For example, if $\lambda\in H^1(X,\Omega^{[1,2>}_X)^*$, thus
$\lambda$ is  a linear function, and
$\chi(z)=\sum_n\chi_{n}z^{-n-1}$, then
\[
\chi(\lambda)(z)=\sum_{n\in ZZ}\lambda(\chi_n)z^{-n-1}\text{ or
}\chi(\lambda)_{(n)}=\lambda(\chi_n).
\]

Denote by $Mod_{\chi}-\cD^{ch,tw}_X$ the full subcategory of
$Mod-\cD^{ch,tw}_X$ consisting of those $\cD^{ch,tw}_X$-modules
where $H_X$ acts according to the character $\chi$.

We will say that a character $\chi(z)\in
H^1(X,\Omega^{[1,2>}_X)((z))$ has  {\em regular singularity} if
$\chi(z)=\chi_0 z^{-1}+\chi_{-1}+\chi_{(-2)}z+\cdots$.

If $\cM\in Mod-\cD^{ch,tw}_X$, then according to Theorem
\ref{zhutcdo}, $\cM_0$ is a $\cD^{tw}_X$-module (even though $\cM$
is filtered and not graded, the fact that the Zhu algebra acts on
the top filtered component remains obviously true). If, in addition,
$\cM  \in   Mod_{\chi}-\cD^{ch,tw}_X$ for some $\chi(z)$ with
regular singularity, then the action of $\cD^{tw}_X$ factors through
the projection $\cD^{tw}_X\rightarrow \cD^{\chi_0}_X$, see
(\ref{specializ-us-tw}), and we obtain a functor
\begin{equation}
 \label{functzerocomp}
\Phi:\; Mod_{\chi}-\cD^{ch,tw}_X\rightarrow Mod-\cD^{\chi_0}_X,
\end{equation}
where $Mod-\cD^{\chi_0}_X$ stands for the category of
$\cD^{\chi_0}_X$-modules.

The locally trivial version
\begin{equation}
 \label{functzerocomploc}
\Phi:\; Mod_{\chi}-\twcdoloc\rightarrow
Mod-\stackrel{\circ}{\cD}^{\chi_0}_X
\end{equation}
is immediate.

The purpose of this section is to prove the following theorem.
\begin{thm}
 \label{equi-categ-tw}
(1)   The category $Mod_{\chi}-\cD^{ch,tw}_X$ consists of only one
object, $\{0\}$, unless $\chi(z)$ has regular singularity.

(2) If $\chi(z)$ has regular singularity, then  functor
(\ref{functzerocomp}) is exact and establishes   an equivalence of
categories.

(3) Assertions (1,2) remain valid upon replacing $\twcdo$ with
$\twcdoloc$.
\end{thm}

\subsection{Proof of Theorem  \ref{equi-categ-tw}.}

Assertion (1) is obvious for if $\chi(z)$ has an irregular
singularity, then condition (3) of the definition of a
$\cD^{ch,tw}_X$ is violated for the subsheaf $\uH_X$.

In order to prove assertion (2) we will construct the left adjoint
to (\ref{functzerocomp}) and show that it is a quasi-inverse of
(\ref{functzerocomp}).

\subsubsection{The left adjoint to (\ref{functzerocomp})}
We  begin by constructing the left adjoint functor locally.

Denote $Mod_\chi {\rm -}{\cD}^{ch,tw}_X(U)$ the category of filtered
$\cD^{ch,tw}_X(U)$-modules  defined
 by analogy with   $Mod_{\chi}{\rm -}\cD^{ch,tw}_X$.

The functor $M \mapsto M_0$ from  $Mod_{\chi}{\rm{-}}\twcdo(U)$ to
$Mod {\rm{-}}{\cD}^{\chi_0}_X(U)$ admits a left adjoint $Zhu_\chi$.
 It is constructed as follows.

Let $F$ be a $\Gamma(U, \cD^{\chi_0}_X)$-module. In particular, it
is a $\Gamma(U, \cD^{tw}_X)$-module, by pullback; therefore we may
apply the functor $Zhu_V$, see section \ref{zhu-left-adj-section},
to it.
 We define
 \begin{equation}
 \tilde{F} = Zhu_V (F),  \ \ V = \cD^{ch,tw}_X(U).
 \end{equation}

For a graded $\Gamma(U, \twcdo)$-module $N$ denote $K_{\chi}(N)$ to
be the subspace spanned by vectors of the form
$$
(\la^*_k)_{(n)}m -  \la^*_k(\chi_n)m,
$$
 where
  $m\in N$, $n\leq -1$, $1\leq k \leq \dim H^1(X, \Omega_X^{[1,2>})$.
It is easy to see that $K_{\chi}(N)$ is a submodule of $N$.

Define
\begin{equation}
Zhu_\chi(F) = \tilde{F}  /  K_{\chi}(\tilde{F}).
\end{equation}
The character $\chi(z)$ having regular   singularity, conditions
(\ref{generated_by_Mzero},\ref{preserved_by_Mzero}) are satisfied;
by construction of the $Zhu_V$-functor, sect\@.
\ref{zhu-left-adj-section}, Condition (3) of a $\Gamma(U,
\twcdo)$-module is satisfied.

 Any $\cD^{\chi_0}(U)$-module map $f: F\to F'$ extends uniquely to a map $Zhu(f):  Zhu_\chi(F)  \to Zhu_\chi(F')$.
Therefore, the functor $Zhu_\chi$ is the left adjoint to the functor
$ M \to M_0$.

Now we proceed to define a sheaf version of $Zhu_\chi$.

If $\cM$ is a $\cD^{\chi_0}_X$-module, let us denote by $\cZ
hu_{\chi}(\cM)$ the sheaf associated to the presheaf
\begin{equation}
U \mapsto Zhu_{\chi, U} \cM(U)
\end{equation}
with restriction maps extended uniquely from those of $\cM$.  It is
clear that $\cZ hu_{\chi}(\cM)$ is a sheaf of
 $\cD^{ch, tw}_X$-modules.
Since maps extend uniquely, this extends to a functor
\begin{equation}
\cZ hu_{\chi}: Mod\rm{-}\cD_X^{\chi_0} \to
Mod_{\chi}\rm{-}\cD_X^{ch, tw}
\end{equation}
left adjoint to the functor \ref{functzerocomp}.

\subsubsection{The quasi-inverse property} We have to show the
following two functor isomorphisms
\begin{equation}
\Phi\circ \cZ hu_{\chi}\stackrel{\sim}{\rightarrow}
\text{Id}_{Mod\rm{-}\cD_X^{\chi_0}}, \label{ovbvios iso-funct}
\end{equation}
\begin{equation}
\cZ hu_{\chi}\circ \Phi\stackrel{\sim}{\rightarrow}
\text{Id}_{Mod_{\chi}\rm{-}\cD_X^{ch, tw}}, \label{notovbvios
iso-funct}
\end{equation}
The first is obviously true, because by construction the functors
are actually equal: $\Phi\circ \cZ hu_{\chi}=
\text{Id}_{Mod\rm{-}\cD_X^{\chi_0}}$. Let us now prove
(\ref{notovbvios iso-funct}).

\bigskip

Let  $U\subset X$ be a suitable for chiralization open subset of
$X$, $A = \Gamma(U, \cO_X)$,
 $\{ \de^i \} $ be an abelian basis for $A$-module $\Gamma(U,
 \cT_X)$,
and $\{  \omega_i \}$   the dual basis of $\Gamma(U, \Omega^1_X)$.

Let $V= \Gamma(U, \twcdo)$, $M$ a $V$-module

Fix a splitting $s: \cT \to V_1$. We will identify $\cT_U$ and
$s(\cT_U) \subset V_{\leq 1}$.

We will denote the  $k$th  mode of $s(\de^i)$ (resp. $\omega_i$) by
$\de_{i, k}$ (resp. $\omega_{i, k}$.)

Let   $P$ denote the polynomial algebra in variables  $\set{D^i_{-
n}, \Omega^i_{-n}, \ n>0, \ 1\leq i\leq \dim X}$

Define the map $a: P \to \End M$, $a(D^i_{- n}) = \de_{i,-n}$,
$a(\Omega^i_{-n})=\omega_{i, -n} $.

Choose any total  order $\succeq$ on the set of variables that
satisfies $D^i_{- n} \succeq \Omega^j_{ -m}\succeq 1$ for all $m>0$,
$n>0$, $i,j$, and $A_{-n}\succeq B_{-m}$ if $n>m$ for $A$ and $B$
being either $D^i$ or $\Omega^j$.

Define $\Psi: P\tensor M_0   \to  M $ as follows:
\begin{eqnarray}
\label{module_identification_map}
  x^1  x^2 \dots x^k \tensor m &\mapsto & a(x^1)  a(x^2) \dots a(x^k) m
 \end{eqnarray}
 where $x^1\succeq x^2\succeq \dots \succeq x^k$; for $k=0$ set $\Psi$ to be the identity map of $M_0$.

\begin{lem} \label{mod-str}
 Suppose $M$ is a filtered $ \Gamma(U, \twcdo)$-module generated
by $M_0$, on which $H_X$
 acts via the character $\chi(z)\in H^1(X, \Omega_X^{[1,2>}[[z]]z^{-1}$
Then:

\indent (A) The map (\ref{module_identification_map}) is a vector
space isomorphism.

\indent (B) If $N\subset M$ is a non-zero submodule, then $N\cap
M_0$ is also non-zero.
\end{lem}

\bigskip

{\bf Proof of Lemma \ref{mod-str}.} (A)  Map
(\ref{module_identification_map}) is surjective by the assumption.
% construction of the left adjoint functor, see  sect.\ \ref{zhu-left-adj-section}.
 To prove
injectivity, extend $\succeq$ to the lexicographic order on the set
of monomials $x^1 x^2 \dots x^k \tensor m$. Let
$\gamma\in\text{Ker}\Psi$ and $\gamma=\gamma_{0}+\cdots$, where
$\gamma_0$ is the leading (w.r.t. the lexicographic ordering)
non-zero term.  Write $\gamma_0=x^1 x^2 \dots x^k \tensor m$. Then
$$
y^1y^2\cdots y^k\Psi(\gamma)=0,
$$
where we choose $y^j$ to be $\frac{1}{n}\partial_{i,n}$ if $x^j=\Omega^{i}_{-n}$
or $\frac{1}{n}\omega_{i,n}$ if $x^j=D^i_{-n}$. The relations of Lemma
\ref{local_model_chir} imply that
$[\partial_{i,n},\omega_{j,-m}]=n\delta_{ij}\delta_{nm}$, and so
thanks to (\ref{generated_by_Mzero})
$$
y^1y^2\cdots y^k\Psi(\gamma)=\frac{\partial}{\partial
x^1}\frac{\partial}{\partial x^2}\cdots\frac{\partial}{\partial
x^k}(x^1x^2\cdots x^k)\Psi(m).
$$
Therefore $\Psi(m)=0$, but the restriction of $\Psi$ to $M_0$ being
the identity, $m$ has to be zero, hence $\gamma=0$, as desired.

Proof of item (B) is very similar: one has to pick a non-zero
$\gamma\in N$ of the lowest degree, and then apply to the highest
degree term, $\gamma_0$, and appropriate $y^1y^2\cdots y^k$ so as to
produce a non-zero element of $N\cap M_0$. $\qed$

\bigskip

Theorem \ref{equi-categ-tw}  follows from Lemma \ref{mod-str}
easily. We have the adjunction morphism
\begin{equation}
\cZ hu_{\chi}\circ \Phi\rightarrow
\text{Id}_{Mod_{\chi}\rm{-}\cD_X^{ch, tw}}, \label{functor-morphism}
\end{equation}
hence
\begin{equation}
\label{transf-of-ident} \cZ hu_{\chi}\circ \Phi(\cM)\rightarrow\cM.
\end{equation}
for each $\cD^{ch,tw}_X$-module $\cM$. The restriction of
(\ref{transf-of-ident}) to $\cM_0$ is the identity.
 By construction, sect.\ \ref{zhu-left-adj-section},
 $ \cZ hu_{\chi}\circ \Phi(\cM)$
 is generated by $\cM_0 = \Phi(\cM)$.
Therefore, due
 to
Lemma  \ref{mod-str}, it is an isomorphism, hence
(\ref{functor-morphism}) is a functor isomorphism. This proves
(\ref{notovbvios iso-funct}).

 $\Phi$ is exact because it is an equivalence of categories; alternatively,
 the exactness follows, immediately, from
Lemma~\ref{mod-str}. The proof of Theorem \ref{equi-categ-tw} (1, 2)
is completed. The locally trivial case, i.e., assertion (3) is
proved in the same way. $\qed$

\begin{rem}
The condition that each $\cM$ be generated by $\cM_0\subset\cM$ in
the definition of a $\twcdo$-module looks unnecessarily restrictive.
Indeed, one can do without it at least when $\twcdo$ is `nice.'

There is an obvious version of the definition of a $\twcdo$-module,
where the generation by $\cM_0\subset\cM$ is replaced with the
requirement that filtration (\ref{filtr-module-2}) exist. Call a
$\twcdo$ {\it locally trivial} if locally on $X$ there is an abelian
basis  $\tau^{(1)},\tau^{(2)},...\subset\cT_X$ and its lift to
$\hat{\tau}^{(1)},\hat{\tau}^{(2)},...\subset(\twcdo)_1$ so that
$\hat{\tau}^{(i)}_{(n)}\hat{\tau}^{(j)}=0$ for all $i,j$ and $n\geq
0$. One can show the following version of Theorem
\ref{equi-categ-tw} is valid for a locally trivial $\twcdo$:

the functor
\begin{align*}
{\rm Sing}:\; Mod_{\chi}-\twcdo&\rightarrow
Mod-\cD^{\chi}_X\\
\cM&\mapsto {\rm Sing}\cM\stackrel{{\rm def}}{=}\{m\in\cM:\,
a_{n}m=0\text{ for all }a\in\twcdo, n>0\}
\end{align*}
is an equivalence of categories. We are planning to return to this
topic in a susequent publication.
\end{rem}

\section{Example: chiral modules over flag manifolds}
\label{Example: chiral modules over flag manifolds}

\subsection{ Sheaf cohomology realization of various $\ghat$-modules. }
\label{results-on-sheaf-coho}

\subsubsection{Bernstein-Beilinson localization. } \label{review-of-bern-beil-loc}
Let $G$ be a complex simple Lie group, $B,B_{-}\subset G$ a generic pair of  Borel subgroups,
$\fg=\text{Lie}\,G$, and $X=G/B_{-}$, the flag manifold. Consider the Beilinson-Bernstein
\cite{BB1} localization functor
\begin{equation}
\label{beil-lern-loc} \Delta: Mod_{ch(\lambda)}-\fg\rightarrow Mod-\cD^{\lambda}_X,
\end{equation}
where we regard $\lambda\in H^{1}(X,\Omega^{1}_{X}\rightarrow\Omega^{2,cl}_X)$ as a weight,
i.e., an element of the dual to a Cartan subalgebra of $\fg$, cf. sect.
\ref{Example:_flag_manifolds.}, especially (\ref{on-flag-trivial}), and $Mod_{ch(\lambda)}-\fg$
is the full subcategory of the category of $\fg$-modules with central character $ch(\lambda)$;
the latter is determined naturally by $\lambda$ and assigns to a central element the number by
which it acts on a module with highest weight $\lambda$. Functor (\ref{beil-lern-loc}) is an
equivalence of categories if $\lambda$ is dominant regular \cite{BB2}.

To see some examples, denote by $V_\lambda$ the simple finite
dimensional $\fg$-module with highest weight $\lambda$, $M_\lambda$
the Verma module with highest weight $\lambda$,  $M_{\lambda}^c$ the
corresponding contragredient Verma module.  We have
\begin{equation}
\label{localiz-of-fin-dim} \Delta(V_\lambda)=\cO(\lambda),
\end{equation}
\begin{equation}
\label{localiz-of-contr-verma}
\Delta(M^c_\lambda) =  i_* i^*\cO(\lambda),
\end{equation}
where $\cO(\lambda)$ is the sheaf of sections of the line bundle
$G\times_{B_-}\Cplx$,
$X_e\stackrel{\text{def}}{=}\bar{B}\subset X$ is the big cell,
$i: X_e \hookrightarrow X$.

\subsubsection{Chiralization.} \label{chiraliz}

Recall that each TCDO on $X$ is locally trivial, see (\ref{on-flag-trivial}). Having fixed
$\chi=\chi(z)\in \fh((z))$ with regular singularity, we obtain the functor
\begin{equation}
\label{g-mod-dchir-mod} \cZ hu_{\chi}\circ\Delta: Mod_{ch(\chi_0)}-\fg\rightarrow
Mod_{\chi}-\cD^{ch,tw}_X,\chi_0=\text{res}_{z=0}\chi(z),
\end{equation}
which is an equivalence of categories if $\chi_0$ is dominant regular (\cite{BB1,BB2} and
Theorem~\ref{equi-categ-tw}.)

According to Lemma \ref{tw-verte-g-f}, there is a vertex algebra
morphism
\begin{equation}
\label{vertex-alg--wakimoto-again} \rho:V_{-h^{\vee}}(\fg)\rightarrow \Gamma(X, \cD^{ch,tw}_X)
\text{ s.t. }\rho(\fz(V_{-h^{\vee}}(\fg))\subset H_X.
\end{equation}

Hence $\cZ hu_{\chi}\circ\Delta(M)$
 is a sheaf of
$V_{-h^{\vee}}(\fg)$-modules with central character $\chi\circ\rho$, where $\chi$ is understood
as in (\ref{def-of-char}). In particular, $\Gamma(X_e,\cZ
hu_{\chi}\circ\Delta(M^{c}_{\chi_0}))$ is a {\it Wakimoto module  of critical level}
\cite{W,FF1,F1}. Indeed, according to (\ref{localiz-of-contr-verma}),
$\Gamma(X_e,\Delta(M^{c}_{\chi_0}))$ is but the space of functions on the big cell $X_e$
carrying an action of $\fg$ twisted by $\lambda$; the definition of the functor $\cZ hu$ in
these circumstances simply mimics the Feigin-Frenkel definition of the  Wakimoto module of
critical level with highest weight $\chi_0$.

 Since
$\Gamma(X,\Delta(M^{c}_{\chi_0}))=M^{c}_{\chi_0}$, we see that the
space of global sections $\Gamma(X,\cZ
hu_{\chi}\circ\Delta(M^{c}_{\chi_0}))$ is the same Wakimoto module
of critical level.

It follows from (\ref{localiz-of-fin-dim},
\ref{localiz-of-contr-verma}) that
$\Gamma(X_{e},\Delta(V_{\chi_0}))=M^{c}_{\chi_0}$, and so
$\Gamma(X_e,\cZ hu_{\chi}\circ\Delta(V_{\chi_0}))$ is also a
Wakimoto module of critical level. What can we say about its space
of global sections?

It is easy to see \cite{MS} that
\begin{equation}
\label{charact-glob-sect}
 \Gamma(X,\cM) \text{ is the maximal
}\fg-\text{integrable submodule of }\Gamma(X_e,\cM)
\end{equation}
Conjecturally, the maximal $\fg$-integrable submodule of a Wakimoto
module of critical level -- and arbitrary highest weight -- is an
irreducible $\widehat{\fg}$-module; therefore $\Gamma(X,\cZ
hu_{\chi}\circ\Delta(V_{\chi_0}))$ is also expected to be
$\widehat{\fg}$-irreducible. We will see how this comes about in the
case where either $\chi_0$ is a regular dominant highest weight --
as has been  assumed so far -- or $\fg=sl_2$ and $\chi(z)=\chi_0/z$,
where $\chi_0$ is an arbitrary integer.

Continuing under the assumption that $\chi_0$ is a regular dominant
integral highest weight we obtain a map
\begin{equation}
\label{Vchi-to-chiral-1} \Delta(V_{\chi_0})\rightarrow\cZ
hu_{\chi}\circ\Delta(V_{\chi_0}),
\end{equation}
hence a map
\begin{equation}
\label{Vchi-to-chiral-2} V_{\chi_0}\rightarrow\Gamma(X,\cZ
hu_{\chi}\circ\Delta(V_{\chi_0})),
\end{equation}

Introduce  the {\it Weyl module of critical level}
$\VV_{\lambda}=\text{Ind}^{\widehat{\fg}}_{\widehat{\fg}_{\leq}}V_{\lambda}$,
where $\fg_\leq$operates on $V_\lambda$ via the evaluation map
$\fg_\leq\rightarrow\fg$, and $K\mapsto -h^{\vee}$, cf.
sect.~\ref{Affine vertex algebras.}. Note that $\VV_0$ is nothing
but the vertex algebra $V_{-h^{\vee}}(\fg)$.

The universality property of induced modules implies that
(\ref{Vchi-to-chiral-2}) uniquely extends to a $\ghat$-morphism
\begin{equation}
\label{Vchi-weyl-to-chiral-2} \VV_{\chi_0}\rightarrow\Gamma(X,\cZ
hu_{\chi}\circ\Delta(V_{\chi_0})).
\end{equation}

This map has kernel, because $\VV_{\lambda}$ carries an action of
the center, $\fz(V_{-h\check{ }}(\fg))$, see Lemma
\ref{tw-verte-g-f}. Define the {\it restricted Weyl module of
central character $\chi(z)$} to be
\begin{equation}
\label{def-of-retrs-weyl-mod}
\VV_{\chi(z)}=\VV_{\chi_0}/\{(p_{(n)}-\chi(\rho(p))_{(n)})v,\;p\in\fz(V_{-h\check{
}}(\fg)),v\in\VV_{\chi_0}.\}
\end{equation}
Then (\ref{Vchi-weyl-to-chiral-2}) factors through
\begin{equation}
\label{Vchi-weyl-restr-to-chiral-2}
\VV_{\chi(z)}\rightarrow\Gamma(X,\cZ
hu_{\chi}\circ\Delta(V_{\chi_0})).
\end{equation}

Frenkel and Gaitsgory \cite{FG3} have proved that $\VV_{\chi(z)}$ is an irreducible
$\ghat$-module.

\begin{thm}\label{descr-of-glob} If $\chi_0$ is regular dominant,
then  map (\ref{Vchi-weyl-restr-to-chiral-2}) is an isomorphism. In
particular,
 $\Gamma(X,\cZ
hu_{\chi}\circ\Delta(V_{\chi_0}))$ is an irreducible $\ghat$-module.
\end{thm}

\bigskip

 Before we continue, let us note that for any smooth variety $X$, even though
$\stackrel{\circ}{\cD}^{ch,tw}_X$ is graded, objects of
$Mod_{\chi}-\twcdoloc$ tend to be only filtered, because quotienting
out by the character $\chi(z)$ does not respect the grading --
except when
\begin{equation}
\label{grad_cond} \chi(z)=\frac{\chi_0}{z},\; \chi_0\in
H^{1}(X,\Omega_X^{1,cl}).
\end{equation}

If $\chi_0$ is integral and $\cL$ is the invertible sheaf of
$\cO_X$-modules with Chern class represented by $\chi_0$, then
Theorem \ref{equi-categ-tw} reads: $Mod_{\chi}-\twcdoloc$ is
equivalent to $Mod-\cD^{\cL}_X$, where $\cD^{\cL}_{X}$ is the
algebra of differential operators acting on $\cL$.  In particular,
associated to $\cL\in Mod-\cD^{\cL}_X$ is $\cZ hu_{\chi}(\cL)\in
Mod_{\chi}-\twcdoloc$. The grading of $\twcdoloc$ induces that of
$\cZ hu_{\chi}(\cL)$:
$$
\cZ hu_{\chi}(\cL)=\cZ hu_{\chi}(\cL)_0\oplus \cZ
hu_{\chi}(\cL)_1\oplus\cdots \text{ where }\cZ hu_{\chi}(\cL)_0=\cL.
$$
Denote
$$
\cL^{ch}=\cZ hu_{\chi}(\cL)
$$
and think of it as  {\it chiralization} of $\cL$.

Suppose now $\fg=sl_2$; then $G/B=\prline$, $H^{1}(X,\Omega_X^{1,cl})=\Cplx$, and we let
$\chi(z)=n/z$, $n\in\ZZ$. We have $\Delta(V_n)=\cO(n)$ if $n\geq 0$ and, independently of the
sign of $n$, $\cO(n)$ is a $\cD_{\prline}^{n}$-module. Therefore, in accordance with the remark
above, we denote by $\cO(n)^{ch}$ the sheaf $\cZ hu_{n/z}(\cO(n))$. The sheaf $\cO(0)^{ch}$ was
one of the first examples of a CDO, and it appeared in \cite{MSV} under the name of the {\it
chiral structure sheaf}.

In this situation, Theorem \ref{descr-of-glob} specializes and
extends as follows:

\begin{thm}
\label{descrforsl2} Let $\LL_n$ be the unique irreducible highest
weight module over $\widehat{sl}_2$ at the critical level with
highest weight $n$. Then

(i) If $n\in\{0,1,2,...\}$, then there are $\widehat{sl_2}$-module
isomorphisms
$$
H^0(\prline,\cO(n)^{ch})\stackrel{\sim}{\rightarrow}H^1(\prline,\cO(n)^{ch})\stackrel{\sim}{\rightarrow}\LL_n.
$$
(ii) If $n\in\{-2,-3,-4,...\}$, then there are
$\widehat{sl_2}$-module isomorphisms
$$
H^0(\prline,\cO(n)^{ch})\stackrel{\sim}{\rightarrow}H^1(\prline,\cO(n)^{ch})\stackrel{\sim}{\rightarrow}\LL_{-n-2}.
$$
(iii) If $n=-1$, then
$$
H^0(\prline,\cO(-1)^{ch})=H^1(\prline,\cO(-1)^{ch})=0.
$$
\end{thm}

\bigskip
\subsection{Proofs}\label{proofofresonsheafcoho} \nopagebreak
\subsubsection{Proof of Theorem \ref{descr-of-glob}.}
Our discussion will  heavily rely on results of \cite{FG1,FG2,FG3}.
Denote by $\cO^{\text{crit}}$ the version of the $\cO$-category of
$\ghat$ modules at the critical level, where all modules $A$ are
assumed to be filtered, $A=\cup_{i=-\infty}^{+\infty}F_iA$, in such
a way that $F_iA=0$ if $i<<0$ and $\fg\otimes t^j(F_iA)\subset
F_{i-j}A$.

Denote by $\fz$ the center of the completed universal enveloping algebra of $\widehat{\fg}$ at
the critical level \cite{FF2}. Any object of $\cO^{\text{crit}}$ is a $\fz$-module. Denote by
$\cO^{\text{crit}}_{\lambda(z)}$ the full subcategory of $\cO^{\text{crit}}$ where $\fz$ acts
according to the character $\lambda(z)$.

Finally, let $\cO^{\text{crit},G}_{\lambda(z)}$ be the full
subcategory of $\cO^{\text{crit}}_{\lambda(z)}$ consisting of
$\fg$-integrable modules. Note that,  by definition,
$\VV_{\chi(z)}\in\text{Ob } \cO^{\text{crit},G}_{\lambda(z)}$
provided $\chi(z)$ and $\lambda(z)$ match, i.e.,
$\lambda=\chi\circ\rho$, see (\ref{vertex-alg--wakimoto-again}).
Likewise, $\Gamma(X,\cZ hu_{\chi}\circ\Delta(V_{\chi_0}))\in\text{Ob
} \cO^{\text{crit},G}_{\lambda(z)}$ thanks to
(\ref{charact-glob-sect}). It is a fundamental result of Frenkel and
Gaitsgory \cite{FG3} that $\cO^{\text{crit},G}_{\lambda(z)}$ is
semi-simple, and $\VV_{\chi(z)}$ is its unique irreducible object.
This implies that map (\ref{Vchi-weyl-restr-to-chiral-2}) is
injective, and it remains to prove surjectivity.

An embedding $\VV_{\chi(z)}\rightarrow A$,
$A\in\text{Ob\;}\cO^{\text{crit},G}_{\lambda(z)}$ is determined by a
{\it singular vector of weight} $\chi_0$, i.e., $v\in A$ such that
(1) $v$ is annihilated by $\fg[t]t\oplus\fn_{+}$, and (2)
$\fh\subset\fg$ operates on $\Cplx v$ according to $\chi_0$. On the
other hand, the semi-simplicity of
$\cO^{\text{crit},G}_{\lambda(z)}$ implies that $\Gamma(X,\cZ
hu_{\chi}\circ\Delta(V_{\chi_0}))$ is a direct sum of copies of
$\VV_{\chi(z)}$. Hence it remains to show that there is a unique up
to proportionality singular vector of weight $\chi_0$ in
$\Gamma(X,\cZ hu_{\chi}\circ\Delta(V_{\chi_0}))$. In fact, more is
true: the entire $\Gamma(X_e,\cZ hu_{\chi}\circ\Delta(V_{\chi_0}))$
contains only one up to proportionality singular vector of weight
$\chi_0$.

To see this, recall that the Wakimoto module  $\Gamma(X_e,\cZ
hu_{\chi}\circ\Delta(V_{\chi_0}))$ is free over $\fn_+[t^{-1}]t^{-1}$ and co-free over
$\fn_+[t]$ with one generator; this fact has been the cornerstone of the Wakimoto module theory
since its inception in \cite{FF1}. A little more precisely, if we let $1$ be the function equal
to 1 on the big cell $X_e\subset X$, then $U(\fn_+[t^{-1}]t^{-1})1\subset \Gamma(X_e,\cZ
hu_{\chi}\circ\Delta(V_{\chi_0}))$ is free and there is  an $\fn_+[t,t^{-1}]$ module
isomorphism
$$
\Gamma(X_e,\cZ hu_{\chi}\circ\Delta(V_{\chi_0}))
\stackrel{\sim}{\rightarrow}
Hom_{\fn_+[t^{-1}]t^{-1}}(U(\fn_+[t,t^{-1}]),U(\fn_+[t^{-1}]t^{-1})1),
$$
where $Hom$ is meant to be the {\it restricted Hom}, which is
defined to be a direct sum of weight components
$$
\oplus_{\alpha,\beta}Hom_{\fn_+[t^{-1}]t^{-1}}(U(\fn_+[t,t^{-1}])_{\alpha},U(\fn_+[t^{-1}]t^{-1})_{\beta}1),
$$
$\alpha$ and $\beta$ varying over the semi-lattice spanned by
positive roots of $\fg$.

It follows that for any $x\in \Gamma(X_e,\cZ
hu_{\chi}\circ\Delta(V_{\chi_0}))$ there is a $u\in U(\fn_+[t])$ so
that $0\neq ux\in U(\fn_+[t^{-1}]t^{-1})1$. Therefore, singular vectors
may occur only in $U(\fn_+[t^{-1}]t^{-1})1$. Weight space
decomposition of the latter is given by
\begin{align}
&U(\fn_+[t^{-1}]t^{-1})1=\oplus_{\alpha}(U(\fn_+[t^{-1}]t^{-1})1)_{\chi_0+\alpha},\nonumber\\
&(U(\fn_+[t^{-1}]t^{-1})1)_{\chi_0+\alpha}=U(\fn_+[t^{-1}]t^{-1})_{\alpha}1.\nonumber
\end{align}
Therefore, $(U(\fn_+[t^{-1}]t^{-1})1)_{\chi_0}$ is one-dimensional
and spanned by $1$, a unique up to proportionality singular vector
of weight $\chi_0$. $\qed$

\bigskip
\subsubsection{ Proof of Theorem \ref{descrforsl2}  }
\label{caseofsl2} Of course item (i) is a particular case of Theorem \ref{descr-of-glob}, but
items (ii, iii) are not. For the reader's convenience we will give an independent proof of all
three items based on representation theory of $\widehat{sl}_2$ as developed in \cite{M}, where
information more complete than in the general case is available; an alternative approach would
be to use \cite{FF1}.

Let $\MM_{\nu}$ be the Verma module over $\widehat{sl}_2$ at the
critical level; this means that $\MM_{\nu}$ is a universal highest
weight  module, where the highest weight vector $v$ satisfies
$h_0v=\nu v; K v=-2v$, cf. sect.\ \ref{Affine vertex algebras.}; we
will also be using some explicit formulas from sect.
\ref{Example:_flag_manifolds.}.

$\MM_{\nu}$ has a unique non-trivial maximal submodule; denote by
$\LL_{\nu}$ the corresponding irreducible quotient.

The Verma module $\MM_{\nu}$ is always reducible, because the
Sugawara operators, which in the vertex algebra notation become
$T_n=(e_{-1}f+f_{-1}e+1/2h_{-1}h)_n$, commute with the action of
$\widehat{sl}_2$. Define the quotient
$$
\MM_{\nu/z}=M_\nu/\sum_{n>0}T_{-n}(\MM_{\nu}).
$$
The module $\MM_{\nu/z}$ is irreducible unless $\nu\in\ZZ-\{-1\}$.
If $\nu=n\in\ZZ-\{-1\}$, then $\MM_{\nu/z}$ is reducible and
contains a unique non-trivial submodule isomorphic to $\LL_{-n-2}$.
We obtain the following exact sequence
\begin{equation}
\label{exactseq} 0\rightarrow \LL_{-n-2}\rightarrow
\MM_{n/z}\rightarrow \LL_n\rightarrow 0.
\end{equation}
The difference between $n$ positive and negative lies in that if
$n\geq 0$, then $\LL_{-n-2}$ is generated, as a submodule, by
$f_0^{n+1}$ applied to the highest weight vector of $\MM_{n/z}$;
therefore, $\LL_n$ is $sl_2$-integrable. On the other hand, if
$n<-1$, then $\LL_{-n-2}$ is generated by $e_{-1}^{-n-1}$ applied to
the highest weight vector of $\MM_{n/z}$; therefore, $\LL_n$ is not
$sl_2$-integrable, but then $\LL_{-n-2}$ is.

 Let us now prove the assertions about the space of global sections
 in (i, ii,iii). In order to compute
$H^0(\prline,\cO(n)^{ch})$, we observe that there is a map
\begin{equation}
\label{Verma-Wakimoto} \MM_n\rightarrow\Gamma(\pole_0,\cO(n)^{ch})
\end{equation}
that sends the highest weight vector  $v\in\MM_n$ to
$1\in\Gamma(\pole_0,\cO(n)^{ch})$, also a highest weight vector.

If $n\geq 0$, $1\in\Gamma(\pole_0,\cO(n)^{ch})$ is annihilated by
$f_0^{n+1}$ (because $1\in\Gamma(\prline,\cO(n))$, and
$\Gamma(\prline,\cO(n))$  is the $(n+1)$-dimensional irreducible
$sl_2$-module.) Therefore, (\ref{Verma-Wakimoto}) factors through
the map
\begin{equation}
\label{irrep-Wakimoto-1}
 \LL_n\rightarrow\Gamma(\pole_0,\cO(n)^{ch}).
\end{equation}
Since $\Gamma(\pole_0,\cO(n)^{ch})$ has the same character as
$\MM_n$, this implies that $\Gamma(\pole_0,\cO(n)^{ch})$ fits into
the following exact sequence
\begin{equation}
\label{exactseq-glo-1} 0\rightarrow \LL_{n}\rightarrow
\Gamma(\pole_0,\cO(n)^{ch})\rightarrow \LL_{-n-2}\rightarrow 0.
\end{equation}
Therefore $\LL_n$ is its unique non-trivial, hence maximal,
submodule, which is $sl_2$-integrable, as it was explained above.
Now (\ref{charact-glob-sect}) implies an isomorphism \linebreak[4]
$H^0(\prline,\cO(n)^{ch})\stackrel{\sim}{\rightarrow}\LL_n$.

In the case where $n<-1$, map (\ref{Verma-Wakimoto}) is an
isomorphism, because the unique non-trivial submodule $\MM_n$ is
generated by $e_{-1}^{-n-1}v$, and map (\ref{Verma-Wakimoto}) sends
the latter to $(\partial_x)_{-1}^{-n-1}1\neq 0$, as formula
(\ref{tw-verte-formulas-sl2-p1}) implies. The brief discussion after
(\ref{exactseq}) shows that if $n<-1$, then the maximal integrable
submodule is $\LL_{-n-2}$ and so is the space of global sections.

Finally, $\Gamma(\pole_0,\cO(-1)^{ch})$ is irreducible and not
integrable, and so the space of global sections is zero.

It remains to show that in each of three cases
$H^1(\prline,\cO(n)^{ch})$ is isomorphic to
$H^0(\prline,\cO(n)^{ch})$. We will achieve that by computing the
Euler character of $\cO(n)^{ch}$ in two different ways.

Since $\cO(n)^{ch}=\oplus_{j\geq 0}\cO(n)_j^{ch}$, we can introduce
the Euler characteristic
$\text{Eu}(\cO(n)_j^{ch})=\dim H^0(\prline,\cO(n)_j^{ch})-\dim H^1(\prline,\cO(n)_j^{ch})$
and the Euler character
$$
\text{Eu}(\cO(n)^{ch})(q)=\sum_{j=0}^{\infty}q^j\text{Eu}(\cO(n)_j^{ch}).
$$
On the other hand, we can similarly consider the formal characters
$ch(H^i(\prline,\cO(n)^{ch}))(q) = \sum_{j\geq 0} q^j \dim
H^i(\prline,\cO(n)^{ch}_j)$, $i=0,1$, and obtain
\begin{equation}
\label{euler-charcter}
\text{Eu}(\cO(n)^{ch})(q)=ch(H^0(\prline,\cO(n)^{ch}))(q)-ch(H^1(\prline,\cO(n)^{ch}))(q).
\end{equation}
The characters of irreducible $\widehat{sl}_2$-modules at the
critical level have been known since [M]; for example,
\begin{equation}
\label{char-ln} ch(H^0(\prline,\cO(n)^{ch}))(q)=ch\,
\LL_n=\frac{n+1}{1-q^{n+1}}\prod_{j=1}^{\infty}(1-q^j)^{-2}\text{ if
}n\geq 0.
\end{equation}

On the other hand, the Euler character $\text{Eu}(\cO(n)^{ch})(q)$
can be computed independently. The sheaf $\cO(n)^{ch}$ carries a
filtration such that the associated graded object is a direct sum of
sheaves $\cO(2s+n)$, $s\in\ZZ$; this is what (\ref{graded_iso_sym})
amounts to in this case. Therefore, we can as well compute the Euler
character of the associated graded object. This is as follows:

\begin{sloppypar}
Informally speaking (cf. \cite{MSV}, sect. 5.8), the local section $(\partial_x)_{-s_1}\cdots
(\partial_x)_{-s_p} x_{-t_1}\cdots x_{-t_q}$ contributes to the graded object the sheaf
$\cO(2p-2q+n)$ sitting in conformal weight $\sum_i(s_j+t_j)$-component. Since
$\text{Eu\;}\cO(2p-2q+n)=2p-2q+n+1$, hence
$\text{Eu\;}\cO(2p-2q+n)+\text{Eu\;}\cO(2q-2p+n)=2(n+1)$, the Euler character of
$\cO(n)^{ch}_j$ equals the number of 2-colored partitions of $j$
 times $(n+1)$. We obtain then
\begin{equation}
\label{formforeul}
\text{Eu}(\cO(n)^{ch})(q)=(n+1)\prod_{j=1}^{\infty}(1-q^j)^{-2}
\end{equation}
\end{sloppypar}
Plugging this and (\ref{char-ln}) in (\ref{euler-charcter}) gives us
\begin{equation}
\label{char-h1} ch(H^1(\prline,\cO(n)^{ch}))(q)=q^{n+1}ch
\LL_n=\frac{(n+1)q^{n+1}}{1-q^{n+1}}\prod_{j=1}^{\infty}(1-q^j)^{-2}\text{
if }n\geq 0.
\end{equation}
We see that $ch(H^1(\prline,\cO(n)^{ch}))(q)$ equals
$ch(H^0(\prline,\cO(n)^{ch}))(q)$ up to an overall power of $q$.
Since an irreducible module is determined by its character, we
conclude that
$H^1(\prline,\cO(n)^{ch})\stackrel{\sim}{\rightarrow}H^0(\prline,\cO(n)^{ch})$,
as desired. Note that the shift by the factor of $q^{n+1}$ means
that $H^1(\prline,\cO(n)^{ch})$ `grows' from the conformal weight
$(n+1)$ component, unlike $H^0(\prline,\cO(n)^{ch})$, which grows
from the conformal weight zero component.

The case of $n<-1$ is handled similarly; an untiring reader will
discover that in this case it is $H^1(\prline,\cO(n)^{ch})$ that
grows from the conformal weight zero component, while
$H^0(\prline,\cO(n)^{ch})$ originates in  the conformal weight
$(-n-1)$ component.

The case where $n=-1$ all the characters in sight are obviously
equal to zero. Theorem \ref{descrforsl2} is now proved. $\qed$

\bigskip

\footnotesize{T.A.: Department of Mathematics, Nara Women's
University, Nara 630-8506, Japan. E-mail address:
arakawa@cc.nara-wu.ac.jp

D.Ch.: Department of Mathematics, University of Southern California,
Los Angeles, CA 90089, USA. E-mail:chebotar@usc.edu

F.M.: Department of Mathematics, University of Southern California,
Los Angeles, CA 90089, USA. E-mail:fmalikov@usc.edu}

\end{document}